\documentclass{article}
\usepackage[utf8]{inputenc}
\usepackage{natbib}
\usepackage{graphicx}
\usepackage{amsmath}
\usepackage{amsfonts}
\usepackage{color}
\usepackage{algorithm}
\usepackage{algpseudocode}
\usepackage{comment}
\algnewcommand\algorithmicswitch{\textbf{switch}}
\algnewcommand\algorithmiccase{\textbf{case}}
\algnewcommand\algorithmicassert{\texttt{assert}}
\algnewcommand\Assert[1]{\State \algorithmicassert(#1)}%
\algdef{SE}[SWITCH]{Switch}{EndSwitch}[1]{\algorithmicswitch\ #1\ \algorithmicdo}{\algorithmicend\ \algorithmicswitch}%
\algdef{SE}[CASE]{Case}{EndCase}[1]{\algorithmiccase\ #1}{\algorithmicend\ \algorithmiccase}%
\algtext*{EndSwitch}%
\algtext*{EndCase}%

\title{Graph Master and Local Area Routes for  Efficient Column Generation for the Capacitated Vehicle Routing Problem with Time Windows}
\author{Udayan Mandal\textsuperscript{ \rm 1,\rm 2}, Amelia Regan\textsuperscript{\rm 2, \rm 3}, Louis Martin Rousseau \textsuperscript{\rm 4}, Julian Yarkony \textsuperscript{\rm 5} \\
\textsuperscript{\rm 1}Stanford University, Palo Alto, California USA\\
\textsuperscript{\rm 2}University of California, Irvine, California USA\\
\textsuperscript{\rm 3}University of Washington, Seattle, Washington USA\\
\textsuperscript{\rm 4} Polytechnique Montreal, Montreal, Quebec, Canada\\
\textsuperscript{\rm 5} Laminaar Optimization Research Group, La Jolla, California USA
}
\date{March 2023}
\begin{document}
\maketitle

\begin{abstract}

 In this research we consider the problem of accelerating the convergence of column generation (CG) for the weighted set cover formulation of the capacitated vehicle routing problem with time windows (CVRPTW). We adapt two new techniques, Local Area (LA) routes and Graph Master (GM) to these problems.  

 LA-routes rely on pre-computing all lowest cost elementary sub-routes, called  LA-arcs, where all customers but the final customer are localized in space. LA-routes are constructed by concatenating LA-arcs where the final customer in a given LA-arc is the first customer in the subsequent LA-arc. To construct the lowest reduced cost elementary route during the pricing step of CG we apply a Decremental State Space Relaxation/time window discretization method over time, remaining demand, and customers visited; where the edges in the associated pricing graph are LA-arcs.

To accelerate the convergence of CG we use an enhanced GM approach.  We map each route generated during pricing to a strict total ordering of all customers, that respects the ordering of customers in the route; and somewhat preserves spatial locality.  Each such strict total ordering is then mapped to a multi-graph where each node is associated with a tuple of customer, capacity remaining, and time remaining. Nodes are connected by feasible LA-arcs when the first/last customers in the LA-arc are less/greater than each intermediate customer (and each other) with respect to the total order. The multi-graph of a given route can express that route and other ``related" routes; and every path from source to sink describes a feasible elementary route. Solving optimization over the restricted master problem over all multi-graphs is done efficiently by constructing the relevant nodes/edges on demand. 
\end{abstract}
\section{Introduction}
In this paper we consider an approach for improving the efficiency of column generation (CG) \citep{barnprice,cuttingstock} methods for solving vehicle routing problems (VRP)\citep{Desrochers1992}.  Our approach can be adapted to other problems solved with CG where pricing is a resource constrained shortest path problem (RCSPP) or an elementary RCSPP (ERCSPP) \citep{irnich2005shortest}. Our paper expands on the recent work on Local Area (LA) routes \citep{mandal2022local_2}, which is applied to the Capacitated Vehicle Routing Problem (CVRP) and can be applied to more elaborate VRP. While we can consider more elaborate problems, we apply our work to CVRPTW for concreteness and ease of exposition. 

A CVRPTW problem instance is described using the following terms: 
\textbf{(a)} a depot located in space; \textbf{(b)} a set of customers located in space, each of which has an integer demand, service time, and time window; and \textbf{(c)} an unlimited number of homogeneous vehicles with integer capacity.  A solution to CVRPTW assigns vehicles to routes, where each route satisfies the following: \textbf{(1)}  The route starts and ends at the depot. \textbf{(2)} The total demand of the customers on the route does not exceed the capacity of the vehicle. \textbf{(3)} The cost of the route is the total distance traveled. \textbf{(4)} The vehicle leaves a customer the same number of times as it arrives at that customer (no teleportation of vehicles). 
\textbf{(5)} Service starts at a given customer within the time window of that customer. \textbf{(6)} No customer is serviced more than once on a route.
The CVRPTW solver selects a set of routes with the goal of minimizing the total distance traveled while ensuring that each customer is serviced at least once and on time. 

Expanded linear programming (LP) relaxations for CVRPTW are celebrated for their tightness relative to compact formulations. In expanded formulations each route is associated with a variable and optimization is treated as weighted set cover.  Since the number of routes can grow rapidly in the number of customers, enumerating all routes is not feasible for large problems.  Hence optimization is done using CG \citep{barnprice,Desrochers1992} which imitates the revised simplex, where the pricing problem is solved by the user and is an ERCSPP \citep{irnich2005shortest}.

Recent work on LA-routes  can be used to accelerate pricing \citep{mandal2022local_2} for CVRP. An LA-route can be a non-elementary route, but must not contain cycles localized in space. Note that a cycle is a section of a route consisting of the same customer at the start and end of the section. Localized cycles in space are cycles consisting of customers that are all spatially close to one another. LA-routes are a tighter version of the celebrated neighborhood-routes (ng-routes) \citep{baldacci2011new}.  LA-routes for CVRP rely on pre-computing the lowest cost elementary sub-route where all customers but the final are localized in space.  LA-routes are constructed by concatenating LA-arcs where the final customer in a given LA-arc is the first customer in the subsequent LA-arc. A Decremental State Space Relaxation  (DSSR) \citep{righini2008new}  is applied to enforce elementarity over LA-routes so as to construct the lowest reduced cost elementary route during the pricing step of CG.

LA-routes also naturally admit a stabilized CG formulation based on Graph Generation (GG) \citep{yarkony2021graph,yarkony2022principled} which we have renamed Graph Master (GM)\footnote{GM alters the CG RMP, not pricing, which the name Graph Generation implies}.  CG stabilization  alters the path of dual solutions generated over the course of CG so as to converge to the optimum more quickly.  In fact the speed up is in the order of hundreds of times faster and increases with problem size.  Specifically GM maps each column produced during pricing to a directed acyclic multi-graph on which every path from source to sink corresponds to a feasible column with identical cost to the actual column. The larger this multi-graph, the more columns can be expressed, but the more difficult computation becomes for the CG restricted master problem (RMP). In the case of CVRPTW, nodes correspond to customer/time/demand remaining and edges are LA-Arcs.  Each multi-graph is associated with a strict total ordering of the customers.  LA-arcs are feasible in a given multi-graph if they satisfy that the first/last customers in the LA-arc are less/greater than the intermediate customers with respect to the ordering. Optimization of this RMP can be done by building up a set of nodes/arcs for the multi-graphs.  Generating new LA-arcs and nodes is a RCSPP that does not have to be elementary and can be solved very efficiently.  This is done by solving pricing (called RMP pricing) over a given multi-graph then adding the nodes/arcs used to the set under  consideration. Since every path in the multi-graph is elementary, this is a RCSPP, which is much easier to solve than its elementary counterpart (ERCSPP) (which still must be solved during pricing but not RMP pricing).  

LA-arcs also permit the efficient inclusion of a slightly weakened form of subset row inequalities (SRI) \citep{jepsen2008subset}.  In their simplest form SRI enforce that the number of selected routes servicing two or more customers in a group of three customers can not exceed one. Such SRI can dramatically tighten the LP relaxation but can increase the time taken by pricing by adding additional resources into the ERCSPP problem. LA-SRI weaken SRI by enforcing that the number of LA-arcs in selected routes including two or more customers in a given group of three customers cannot exceed one. The use of LA-SRI does not alter the difficulty of pricing, since it does not inject additional resources that must be considered.  Since routes in an optimal solution tend not to return to areas where they have already visited, LA-SRI can dramatically tighten the set cover LP relaxation of CVRP.  

Taking advantage of LA-routes and the associated stabilization approach is non-trivial when considering time windows. This is because the optimal ordering of customers in an LA-arc varies with the departure and arrival time from the first/last customers on the LA-arc. We cannot compute the optimal ordering via enumeration for each departure/arrival time as time is highly fine grained.  To circumvent this we introduce a dynamic programming-based approach to generate the efficient frontier trading off latest departure, earliest arrival and distance traveled. This is based on the observation that removing the  first customer from an ordering on an efficient frontier must produce an ordering that lies on the efficient frontier for its respective LA-arc.  Thus we can construct the efficient frontier for all LA-arcs jointly by generating them from smallest to largest (in terms of the number of intermediate customers).  Given an LA-arc and start/end time we can generate the lowest cost arc by considering only the orderings in the efficient frontier.  To permit the solution to the pricing problem efficiently we expand the work of \citep{boland2017continuous} so as to associate the nodes in the CG pricing problem with ranges of time/capacity/feasible customers remaining.  These sets are then divided over the course of pricing until the lowest reduced cost route is elementary and feasible. 

We organize this document as follows.  In Section \ref{Sec_litReview} we review related work.  In Section \ref{Sec_mathFormBasics} we provide the basic description  used for CVRPTW with a focus on the CG solution.  In  Section \ref{sec_graph_master} we provide an accelerated approach to CG for CVRPTW given a mechanism to solve the pricing problem.  In Section \ref{sec_pricing} we describe a fast solution to the pricing problem using LA-arcs.  In Section \ref{sec_la_time} we describe the computation of all optimal orderings for customers in LA-arcs over time using dynamic programming. In Section \ref{sec_exper} we provide experimental validation of our approach.  In Section \ref{sec_conc} we conclude and discuss extensions with an emphasis on branch-cut-price and the integration of LA-SRI \citep{mandal2022local_2}.  
\section{Literature Review}
\label{Sec_litReview}

\subsection{Neighborhood Routes/Decremental State Space Relaxation}

Decremental State Space Relaxation (DSSR) is an efficient mechanism for solving pricing that gradually enforces elementarity as needed, thus lowering the computational burden.  DSSR alternates between \textbf{(1)} generating the lowest reduced cost path partially relaxing elementarity (and enforcing resource feasibility, meaning that the path is feasible with respect to any resource constraints such as time and demand), and \textbf{(2)} augmenting the constraints enforced so as to prevent the current non-elementary solution from being regenerated. Step \textbf{(1)} produces the path with the lowest reduced cost from a super-set of the set of elementary routes; this set decreases in size as DSSR proceeds. Termination of DSSR is achieved when the generated path is elementary, at which point this path is guaranteed to be the lowest reduced cost elementary route. In practice DSSR does not need to generate all such constraints. DSSR encodes constraints by associating each customer $u$ with a set of the other customers called its neighborhood, which is denoted $M_u$. The path generated at a given iteration of DSSR does not include any cycle satisfying the following property: The cycle starts and ends at $u$, and $u \in M_v $ for all intermediate customers $v$. Generating such a path in step \textbf{(1)} is tackled as a dynamic programming problem, which can be alternatively solved using labeling algorithms \citep{Desaulniers2005}.  Given a non-elementary path generated in \textbf{(1)}, a cycle is identified; then in step \textbf{(2)} the neighborhood sets of all intermediate customers are augmented to include the starting/ending customer of the cycle. The solution time of the labeling algorithm can grow exponentially as a function of the maximum size of any neighborhood. 

The Neighborhood route (ng-route) \citep{baldacci2011new} relaxations can be understood as DSSR where we do not grow the set on neighbors but instead initialize the neighbors of each given customer to be the set of customers it is spatially nearby. It is observed that the use of ng-routes can vastly improve the performance of CG based solvers for VRP.  
\subsection{General Dual Stabilization }
The number of iterations of CG required to optimally solve the weighted set cover formulation of CVRP, which is also called the master problem (MP), can be dramatically decreased by intelligently altering the sequence of dual solutions generated \citep{du1999stabilized, rousseau2007interior, Pessoa2018Automation} over the course of CG. Such approaches, called dual stabilization, can be written as seeking to maximize the Lagrangian bound at each iteration of CG \citep{geoffrion1974lagrangean}.  The Lagrangian bound is a lower bound on the optimal solution objective  to the MP that can be easily generated at each iteration of CG.  In CVRP problems the Lagrangian bound is the LP value of the RMP plus the reduced cost of the lowest reduced cost column times the number of customers. Observe that when no negative reduced cost columns exist, the Lagrangian bound is simply the LP value of the RMP. The Lagrangian bound is a concave function of the dual variable vector. The current columns in the RMP provide for a good approximation of the Lagrangian bound nearby dual solutions generated thus far but not regarding distant dual solutions.  This motivates the idea of attacking the maximization of the Lagrangian bound in a manner akin to gradient ascent. Specifically we trade off maximizing the objective of the RMP, and having the produced dual solution be close to the dual solution with the greatest Lagrangian bound identified thus far (called the incumbent solution).  

A simple but effective version of this idea is the box-step method of \citep{marsten1975boxstep}, which maximizes the Lagrangian bound at each iteration of CG such that the dual solution does not leave a bounding box around the incumbent solution.  Given the new solution, the lowest reduced cost column is generated, and if the associated Lagrangian bound is greater than that of the incumbent then the incumbent is updated. The alternative stabilization approach of \citep{Pessoa2018Automation} takes the weighted combination of the incumbent solution and the solution to the RMP and performs pricing on that weighted combination. 
\subsection{Time Window Discretization}
Our pricing solver draws on the insights discussed in \citep{boland2017continuous} with respect to continuous-time service network design problems. In this paper the authors iteratively construct a partition of the times by dividing nodes corresponding to ranges of time. Such networks permit a vehicle to leave a node before it arrives. By iteratively solving the mixed integer linear program and splitting nodes used that cause a violation of temporal feasibility, the authors rapidly solve the optimization. This work is adapted in a CG context in the context of automated warehouses in \citep{haghani2021multi}.  
\subsection{Graph Master/Principled Graph Management}
Graph Master (GM) is an approach to stabilizing CG by permitting each column generated to describe many columns without altering the structure of the pricing problem \citep{yarkony2021graph,yarkony2022principled}.  
To apply GM, we must be able to map any given column to a directed acyclic multi-graph for which any path from source to sink describes a feasible column. This structure is easily satisfied for vehicle routing, crew scheduling problems, and other problems where pricing is a resource constrained shortest path problem. Such multi-graphs are then added to the RMP when the corresponding column is generated during pricing. The use of GM does not weaken the linear programming (LP) relaxation being solved. GM permits the RMP to express a much wider set of columns than those generated during pricing, leading to faster convergence relative to standard CG. GM does not change the structure of the CG pricing problem.  The GM RMP has a specific primal block angular structure that can be efficiently exploited to solve the MP.  

GM has two computational bottlenecks. The first is pricing. The structure of the problems solved using GM is identical to that of standard CG. The second bottleneck is the RMP, which is more computationally intensive in GM than in standard CG given the same number of columns generated. By design GM converges in fewer iterations than standard CG, and hence requires fewer calls to pricing. Therefore when the computation time of GM is dominated by pricing, as opposed to solving the RMP, GM converges much faster than standard CG in terms of time. However GM need not converge faster than standard CG when the GM RMP, rather than pricing, dominates computation.  This issue is addressed by Principled Graph Management (PGM) \citep{yarkony2022principled}.  PGM is an algorithm to solve the GM RMP rapidly by exploiting its special structure.  In calls to the RMP, PGM iterates between solving the RMP over a subset of edges and adding new edges. Computing new edges to consider is done by identifying paths in a given multi-graph with negative total cost, which are then added to the RMP. Such paths are generated as a standard shortest path computation not an ERCSPP.  

\subsection{Local Area Routes}
Local area routes are introduced in \citep{mandal2022local_2} so as to improve the efficiency of pricing, the number of iterations of CG and tighten the LP relaxation. 

\textbf{Pricing with Local Area Routes:} In the context of CVRP, LA-routes rely on pre-computing the lowest cost elementary sub-route (called an LA-arc) for each tuple consisting of the following: \textbf{(1)} a (first) customer where the LA-arc begins, \textbf{(2)} a distant customer (from the first) where the LA-arc ends, and \textbf{(3)} a set of intermediate customers near the first customer. LA-routes are constructed by concatenating LA-arcs where the final customer in a given LA-arc is the first customer in the subsequent LA-arc. A Decremental State Space Relaxation method  \citep{righini2008new,righini2009decremental} is applied over LA-routes to construct the lowest reduced cost elementary route during the pricing step of CG. 

\textbf{Local Area Stabilization:}  LA-routes can be efficiently encoded in a GM based solution for efficiently solving the MP. Specifically each column generated during pricing is mapped to an strict total ordering of the customers consistent with that column. An LA-arc is consistent with an ordering if the first/last customer in the arc come before/after all other customers in the LA-arc in the associated ordering respectively. Each such ordering is then mapped to a multi-graph where nodes correspond to (customer/demand) and edges correspond to LA-arcs consistent with that ordering. Hence any path from source to sink on the multi-graph is a feasible elementary route.  The ordering for a column places customers spatially nearby in nearby positions on the ordering so that routes can be generated so as to permit spatially nearby customers to be visited without traveling far away first. We solve the RMP over these graphs, which has special structure allowing for efficient solution.

\textbf{Local Area Subset Row Inequalities:}  LA-route based solvers can be used to efficiently tighten the standard weighted set cover formulation of CVRP using a variant of subset row inequalities (SRI) \citep{jepsen2008subset}, which do not alter the structure of pricing. LA-SRI marginally weaken standard SRI by enforcing the SRI constraint over LA-arcs (with the last customer ignored in each arc to avoid over-counting).  This then permits their efficient use in pricing and in the stabilized RMP.  The inclusion of LA-SRI does not alter the structure of the pricing, permitting many more LA-SRI constraints to be included in the MP and considered during pricing.  Rounded capacity inequalities (RCI) \citep{archetti2011column} can also be used with LA-arcs inside the stabilized formulation.

\section{Capacitated Vehicle Routing Problem with Time Windows}
\label{Sec_mathFormBasics}
In this section we formulate CVRPTW as a weighted set cover problem, as is common in the operations research literature \citep{costa2019}.  We use $N$ to denote the set of customers, which we index by $u$. We use $N^+$ to denote $N$ augmented with the starting/ending depot (which are co-located), and are denoted $-1,-2$ respectively.  
Each customer $u\in N$ has demand $d_u$. We define the demand at the depot as $0$, which we write formally as $d_{-1}=d_{-2}=0$.  
We use $t^+_u$ and $t^-_u$ to denote the earliest/latest times for the service at customer $u$ to begin.  
The service window for the start/end depot is $t_{-1}^{+}=t_{-2}^{+}=t_0$ and $t_{-1}^{-}=t_{-2}^{-}=0$, where $t_0$ is the length of the time horizon of optimization.  Servicing a given customer $u\in N$ takes $t_u^*$ units of time, where $t_{-1}^*=t_{-2}^*=0$.  
For each $u \in N^+,v\in N^+$ we use $c_{uv}=t_{uv}$ to denote the cost and time required to travel from $u$ to $v$ respectively. For simplicity $t_{uv}$ is altered so as to include the service time at customer $u$, which is $t_u^*$.  Thus we set $t_{uv}\leftarrow t_{uv}+t_{u}^*$ then set $c_{uv}\leftarrow t_{uv}$ (which offsets the cost of a solution by a constant).  Thus $t_{u}^*$ is removed from consideration for the rest of the document.

We are given a set of homogeneous vehicles each with capacity $d_0$.  We use $\Omega$ to denote the set of all feasible routes, which we index by $l$.  We describe $l$ using binary term $a_{uvl}$; where $a_{uvl}=1$  if and only if (IFF) $u$ is succeeded by $v$ in route $l$.  We set $a_{ul}=1$ if $u$ is visited in route $l$, and otherwise set $a_{ul}=0$.  We use $N_l$ to denote the set of customers visited by route $l$.  Below, we define $c_l$ to be the cost of route $l$, which is the total distance traveled on route $l$.
\begin{align}
c_l=\sum_{\substack{u\in N^+\\ v \in N^+}}c_{uv}a_{uvl} \quad \forall l \in \Omega
\end{align}
We now formulate CVRPTW as weighted set cover problem over $\Omega$. Here routes in $\Omega$ are sets and we must cover each customer at least once.   We set decision variable $\theta_l=1$ to indicate that route $l$ is selected in our solution and otherwise set $\theta_l=0$.  We use $\Psi(\Omega)$ to refer to the weighted set cover problem over $\Omega$.  For computational reasons later we use $\Psi(\Omega_R)$ to refer to weighted set cover over a subset of $\Omega$ denoted $\Omega_R$.  Here $\Psi(\Omega_R)$ and $\Psi(\Omega)$ are referred to as the restricted master problem (RMP) and the master problem (MP) respectively.   Below we write CVRPTW as a linear program (LP) in primal and dual form.  
%
\begin{subequations}
\label{primal_master}
\begin{align}
\Psi(\Omega_R)=\min_{\theta \geq 0}\sum_{l \in \Omega_R}c_l\theta_l \label{CVRP_obj}\\
\sum_{l \in \Omega_R}a_{ul}\theta_l\geq 1 \quad \forall u \in N \quad [\pi_u] \label{CVRP_cover}
\end{align}
\end{subequations}
\begin{subequations}
\label{dualForm}
    \begin{align}
        \Psi(\Omega_R)=\max_{\pi \geq 0}\sum_{u \in N}\pi_u \label{dualForm_obj}\\
        c_l-\sum_{u \in N}a_{ul}\pi_u \geq 0 \quad \forall l \in \Omega_R \label{dualForm_con} \quad [\theta_l]
    \end{align}
\end{subequations}
\textbf{Description of the Primal LP:  }In \eqref{CVRP_obj} we minimize the total cost of the routes used.  In \eqref{CVRP_cover} we ensure that each customer is serviced at least once (though an optimal solution services each customer exactly once).  We use $[\pi_u]$ to indicate the dual variable associated with a given primal constraint.  

\textbf{Description of the Dual LP:  }In \eqref{dualForm_obj} we maximize the sum of the dual variables $\pi_u$.  In \eqref{dualForm_con} we enforce that each route has non-negative reduced cost.  We use $[\theta_l]$ to indicate the primal variable associated with a given dual constraint.  

When $\Omega_R\leftarrow \Omega$ \eqref{primal_master} and \eqref{dualForm} are challenging to solve optimally with an off the shelf LP solver as the number of primal variables grows rapidly in the number of customers.  We can not enumerate all such variables much less consider them in optimization.  
Column generation (CG) \citep{cuttingstock,Desrochers1992} solves $\Psi(\Omega)$ without having to consider all of $\Omega$. CG constructs a sufficient subset of $\Omega$ denoted $\Omega_R$ s.t. $\Psi(\Omega_R)$ provides an optimal solution to $\Psi(\Omega)$. 
To construct $\Omega_R$, CG iterates between \textbf{(1)} solving $\Psi(\Omega_R)$, and \textbf{(2)} identifying at least one $l \in \Omega$ with negative reduced cost, which are then added to $\Omega_R$.  We can initialize $\Omega_R$ with one route for each $u\in N$ that services just that single customer.  In many CG applications only the lowest reduced cost primal variable in $ \Omega$ (also referred to as a column) is generated. We write the selection of this column as optimization below using $\bar{c}_l$ to denote the reduced cost of column $l$.
\begin{subequations}
\label{pricing}
\begin{align}
    \min_{l \in \Omega} \bar{c}_l \\
    \bar{c}_l=c_l-\sum_{u \in N}a_{ul}\pi_u  \quad \forall l \in \Omega \label{redCostForm}
\end{align}
\end{subequations}
The operation in \eqref{pricing} is referred to as pricing and solved as an elementary resource constrained shortest path problem (ERCSPP)\citep{irnich2005shortest}, not by enumerating all of $\Omega$, which would be computationally prohibitive. CG terminates when pricing proves no column with negative reduced cost exists in $\Omega$. This certifies that CG has produced the optimal solution to $\Psi(\Omega)$. This ERCSPP must keep track of the following resources \textbf{(1)}  the set of customers already visited; \textbf{(2)} the total demand consumed; \textbf{(3)} the amount of time remaining.

Pricing itself can be written as an mixed integer linear program (MILP) though it is not typically solved as a MILP. We now describe the MILP formulation for pricing for CVRPTW for completeness. We set binary decision variable $x_{uv}=1$ to indicate that the vehicle travels from $u$ directly to $v$, for any $u,v$ pair each in $N^+$.  For each $u \in N$ we use $y_u$ to denote the amount of time remaining when the vehicle departs $u$.  For each $u\in N$ we use $z_{u}$ to denote the amount of capacity remaining immediately prior to starting service at $u$.  Note that $y_u,z_u$ only have values when $u$ is visited in the generated route.  
\begin{subequations}
\begin{align}
\min_{\substack{x \in \{0,1\}\\ y \geq 0\\ z\geq 0}}\sum_{\substack{u\in N^+\\v \in N^+}}(c_{uv}-\pi_u)x_{uv} \label{pricingObj}\\
\sum_{u \in N}x_{-1u}=1 \label{leavOnce_feas}\\
\sum_{v \in N^+}x_{uv}-x_{vu}=0 \quad \forall u \in N \label{flowCon_feas}\\
y_u-t_{uv}+(1-x_{uv})\infty\geq y_v \quad \forall u\in N^+,v\in N^+ \label{timeEdge_feas}\\
z_u-d_u+(1-x_{uv})\infty \geq z_v \label{demEdge_feas} \quad \forall u\in N^+,v \in N^+\\
d_0\geq z_u\geq d_u \label{dem_feas} \quad \forall u \in N^+\\
t^+_u\geq y_u\geq t^-_u \label{time_feas} \quad \forall u \in N^+
\end{align}
\end{subequations}
In \eqref{pricingObj} we seek to minimize the travel distance minus the dual variables of the customers serviced (here $\pi_{-1}$ and $\pi_{-2}$ are defined to be zero).  The objective in \eqref{pricingObj} is the reduced cost of the generated route.  In \eqref{leavOnce_feas} we enforce that the vehicle leaves the starting depot exactly once.  In \eqref{flowCon_feas} we enforce that the vehicle leaves a customer the same number of times as it arrives at that customer.
In \eqref{timeEdge_feas} we enforce that if the vehicle travels from $u$ to $v$ then it departs $v$ no earlier than when it departs $u$ minus the travel time from $u$ to $v$.
In \eqref{demEdge_feas} we enforce that if the vehicle travels from $u$ to $v$ then the amount of capacity remaining prior to servicing $v$ is no more than the amount remaining prior to servicing $u$ minus $d_u$.  In \eqref{dem_feas} we enforce that the vehicle arrives at $u$ with at least enough capacity remaining to service $u$.
In \eqref{time_feas} we enforce time window feasibility.  
%
\section{Graph Master for Efficient Column Generation Optimization} 
\label{sec_graph_master}
In this section we consider a novel Graph Master (GM) \citep{yarkony2021graph,yarkony2022principled} based approach for efficiently solving $\Psi(\Omega)$.  GM is employed so as to limit the number of calls to pricing and therefore the total computation time. 

We organize this section as follows. In Section \ref{laNeighArcs} we describe LA-arcs in the context of CVRPTW.  In Section \ref{sec_super_graph} we produce a multi-graph called the super LA master graph, which is never computed, but can describe all routes needed to solve $\Psi(\Omega)$ optimally.  In Section \ref{sec_graph_master_form} we describe sub-graphs of the super LA master graph used to describe subsets of columns called families.  In Section \ref{sec_solution_appraoch} we use GM to solve the RMP over large families of columns in a computationally efficient manner.
\subsection{Local Area Arcs for CVRPTW}
\label{laNeighArcs}
We now describe LA-arcs in the context of CVRPTW.  For each customer $u\in N$ we use $N_u
\subseteq (N-u)$ to denote the set of customers that are nearest to $u$ (excluding customers that can not be reached from $u$ due to time). We refer to $N_u$ as the set of LA-neighbors of $u$. The starting/ending depot have no LA-neighbors and no customer considers the starting/ending depot to be one of its LA-neighbors.  

We use $P$, which we index by $p$, to denote the set of LA-arcs where $p$ is defined as a tuple $(u_p,v_p,N_p)$.  Each arc consists of the following \textbf{(1)} customer/depot $u_p$; \textbf{(2)} a subset of the set of LA-neighbors of $u_p$, which we denote as $N_p$ (with $N_p^+=(u_p \cup N_p)$); \textbf{(3)} a customer/depot that is not in the set of LA-neighbors of $u_p$ denoted $v_p$. 
We use $a_{up}$ to describe the customers in LA-arc $p$ (excluding the final customer); thus $a_{up}=1$ if $u \in N^+_p$, and otherwise $a_{up}=0$ for all $u\in N$.  We use $d_p$ to describe the total amount of demand serviced on LA-arc $p$ excluding $v_p$ as follows:  $d_p=\sum_{w \in N^+_p}d_w$.  The customers in an LA-arc must be able to be serviced in a single route with regards to vehicle capacity meaning $d_0 \geq d_p$.  For any given $p \in P$ and times $t_1,t_2$ we use $c_{p,t_1,t_2}$ to denote the cost of the lowest cost sequence of customers departing $u_p$ at time no earlier than $t_1$, and departing $v_p$ no later than $t_2$, and visiting $N_p$ as intermediate customers (between $u_p$ and $v_p$).  
Later in the document in Section \ref{sec_la_time} we consider the efficient computation of $c_{p,t_1,t_2}$, however for now we assume that we can compute $c_{p,t_1,t_2}$ efficiently.
\subsection{Super LA Master Graph}
\label{sec_super_graph}
In this sub-section we consider a multi-graph called the super LA master graph. We denote the super LA master graph as $G$, which has node set $I$ and edge set $E$. We describe $I$ as follows. There is one node in $I$ for the source $(-1,d_0,t_0)$, sink $(-2,0,0)$, and each $(u,d,t)$ for which $u \in N,d_0\geq d\geq d_u,t^+_u\geq t\geq t^-_u$.  We describe the edge set $E$ below.  
\begin{itemize}
    \item For each $i=(-1,t_0,d_0)$ and $j=(u,d,t)\in I$  (where $t\leq t_0-t_{uv}$) we connect $i$ to $j$ with an edge of cost $c_{-1,u}$.  Traversing this edge indicates that a vehicle leaves the depot and then arrives at customer $u$  visiting no intermediate customers. The vehicle leaves the depot with $d$ units of capacity remaining and arrives at $u$ with $t$ units of time remaining. Note that this vehicle need not have $d_0$ capacity remaining when leaving the depot, meaning that it is not fully loaded.  
    The customers serviced on this edge are denoted $N_{ijp}=\{ \}$ for $p=\{-1,u,\{ \}\}$ (in other words, no customers are serviced when traveling from the depot to $u$).  
    \item Consider each tuple of $i \in I, j \in I, p\in P$ for which $ i=(u_p,d_1,t_1)$ $j=(v_p,d_2,t_2)$, $c_{p,t_1,t_2}<\infty$, $d_1-d_p\geq d_2$.
    We connect $i$ to $j$ with an edge of cost $c_{p,t_1,t_2}$ where traversing this edge indicates that we leave $u_p$ at time $t_1$, then service $N_p$ in the lowest cost feasible order then go to $v_p$, which we leave at time $t_2$. The customers serviced on this edge are denoted $N_{ijp}=N^+_p$.  Note that $v_p$ is not included in $N_{ijp}$.  Also we define $N_{ijp}$ to be ordered from earliest to latest serviced.
\end{itemize}
Given a path from source to sink  we use  $[e_1,e_2,e_3...|e|]$ to describe the edges in the path.  Here $|e|$ is the total number of edges in the path and $e_k$ describes the $k$'th edge visited and $e_k=(i_k,j_k,p_k)$ are the associated nodes and LA-arc of the $k$'th edge. The customers in the path are visited in the following order
$[N^+_{p_1},N^+_{p_2},N^+_{p_3},N^+_{p_4}...N^+_{p_{|e|}}]$.

Observe that by construction of the graph that the total capacity of customers serviced in this path is no greater than $d_0$ and all time windows are obeyed.  Furthermore if $|N^+_{p_{k1}}\cap N^+_{p_{k2}}|=0$ for all $k_1\neq k_2$ then the path is elementary and hence the path defines an elementary route lying in $\Omega$. Observe that every elementary path from source to sink corresponds to an $l\in \Omega$; and that it will have total cost $c_l$ and service the customers in $N_l$.
Observe that for any $l\in \Omega$ for which, $l$ is the lowest cost route in $\Omega$ servicing $N_l$, then $l$ can be represented in $G$. 
\subsection{Graph Master Structure Exploiting LA-Routes}
\label{sec_graph_master_form}
We consider a directed acyclic multi-graph defined by $l \in \Omega,I^-\subseteq I,P^-\subseteq P$ and vector $\beta^l$.  We define $I^-$ to include the nodes corresponding to the source and the sink in $I$.  
Here $\beta^l$ describes a strict total ordering of the customers/depots.  For a given $u\in N^+$, $\beta^l_u$ indicates the position of $u$ in this ordering. The starting/ending depot have the smallest and largest $\beta^l$ terms respectively.  Given $I^-,P^-,\beta^l$  we describe the edge set as $E^{l,I^-,P^-}$ which is defined as follows. 
\begin{subequations}
\label{def_arcs_family}
\begin{align}
E^{l,I^-,P^-}=\{ (i,j,p) \in E\\
i\in I^-,j \in I^-, p\in P^- \\
\beta^l_{u_p}<\beta^l_w \quad \forall  w \in (v_p\cup N_p) \\ 
\beta^l_{w}<\beta^l_{v_p} \quad \forall  w \in (N_p\cup u_p )\}
\end{align}
\end{subequations}
Observe that since any path from source to sink over $(I^-,E^{l,I^-,P^-})$ is also a path in $G$ then that path satisfies demand and time window feasibility.  Observe that every path from source to sink in the graph $(I^-,E^{l,I^-,P^-})$ describes an elementary route as no customers can be repeated by \eqref{def_arcs_family}.
The set of routes that can be expressed in $(I,E^{l,I,P})$ (made by setting $I^-\leftarrow I $ and $P^-\leftarrow P$) is referred to as $\Omega_l$ where $l \in \Omega_l $ and $\Omega_l \subseteq \Omega$.  We refer to $\Omega_l$ as the family of route $l$.  The set of routes expressible using $(I^-,E^{l,I^-,P^-})$ is denoted $\Omega_{l,I^-,P^-}$, which is a subset of $\Omega_l$.  Observe that longer LA-arcs (longer means larger sets $|N_p|$) can permit more routes to lie in $\Omega_l$.

  Observe that the construction of $\beta^l$ dramatically alters the routes in $\Omega_l$.  The construction of $\beta^l$, which we share with \citep{mandal2022local_2}, is motivated by the observation that customers that are in similar physical locations should be in similar positions on the ordered list.  Having customers in such an order ensures that a route defined over $(I,E^{l,I,P})$ can visit all customers close together in an area without leaving the area, and then coming back. 

We now describe the construction of $\beta^l$ as provided in \citep{mandal2022local_2}.  We initialize the ordering with the members in $N_l\cup -1\cup -2$ sorted in order from first visited to last visited in route $l$. Then, we iterate over $v\in (N-N_l)$, and insert $v$ behind the $u \in (N_l \cup -1)$ for which $c_{uv}$ is minimized \footnote{In experiments for customers closest to the depot we placed them at the end of the list}.  As an aside we note that when the time windows or demand of $u,v$ are such that the route $[-1,u,v,-2]$ is infeasible then $c_{uv}$ is regarded as $\infty$.  Here $[-1,u,v,-2]$ describes the route starting at the starting depot then servicing $u$ then $v$ before going to the ending depot.
\subsection{Graph Master Based Solver}
\label{sec_solution_appraoch}
We now consider our GM based solver for $\Psi(\Omega)$. Our solver consists of two nested loops.  The outer loop solves $\Psi(\Omega)$ by alternating between \textbf{(1)} solving $\Psi(\cup_{l \in \Omega_R}\Omega_l)$ and \textbf{(2)} adding the lowest reduced cost route $l\in \Omega$  to $\Omega_R$.  The solution to $\Psi(\cup_{l \in \Omega_R}\Omega_l)$ is done by the inner loop, which solves $\Psi(\cup_{l \in \Omega_R}\Omega_l)$ by alternating between \textbf{(1)} solving $\Psi(\cup_{l \in \Omega_R}\Omega_{l,I^-,P^-})$  and \textbf{(2)} adding elements to $I^-,P^-$. Specifically, elements are added to $I^-,E^-$ such that the lowest reduced cost route in $\Omega_l$ can be represented for each $l\in \Omega_R$.  The computation of the lowest reduced cost route $l \in \Omega_{\hat{l}}$ for each $\hat{l} \in \Omega_R$ is considered during our general discussion of pricing in Section \ref{sec_pricing}.  We let $P_l,I_l$ denote the LA-arcs and nodes from $P,I$ used in route $l$.  Thus given that $l$ is the lowest reduced cost route in $\Omega_{\hat{l}}$ for some $\hat{l}\in \Omega_R$ we add $I_l,P_l$ to $I^-,P^-$ respectively.  

We describe the efficient solution to $\Psi(\cup_{l \in \Omega_R}\Omega_{l,I^-,P^-})$ which for short hand we refer to as $\Psi(\Omega_R,I^-,P^-)$  using the following decision variables.  For each $l \in \Omega_R, (i,j,p) \in E^{l,I^-,P^-}$ we create one decision variable $x^l_{ijp}$. Here $i,j,p$ are described as follows:  $i=(u,d_1,t_1)$ and $j=(v,d_2,t_2)$ and $p \in P$ with $p=(u,v,N_p)$.
Setting $x^l_{ijp}=1$ indicates that a vehicle leaves $u$ with at least $d_1-d_u$ units of capacity remaining and $t_1$ units of time remaining, then services the customers in $N_p$, then goes to $v$, which it departs at $t_2$ with at least $d_2-d_v$ units of capacity remaining.  Below we formulate $\Psi(\Omega_R,I^-,P^-)$ as a LP with exposition provided below the formulation. 
%
%
\begin{subequations}
\label{my_rmp}
\begin{align}
    \Psi(\Omega_R,I^-,P^-)=\min_{x \geq 0}\sum_{\substack{l \in \Omega_R\\ijp\in E^{lI^-P^-}}}c_{pt_it_j}x^l_{ijp} \label{objTerm}\\
    \sum_{\substack{l \in \Omega_R\\ijp\in E^{lI^-P^-}}}a_{up}x^l_{ijp}\geq 1 \quad \forall u \in N \quad [\pi_u] \label{cover_term}\\
    \sum_{\substack{ijp\in E^{lI^-P^-}}}x^l_{ijp}=\sum_{\substack{jip\in E^{lI^-P^-}}}x^l_{jip} \quad \forall l \in \Omega_R,i \in I^- -((-1,d_0,t_0) \cup (-2,0,0)) \label{flow_term}
\end{align}
\end{subequations}

In \eqref{objTerm} we minimize the total cost of the routes used, by minimizing the total cost of the LA-arcs used.
In \eqref{cover_term} we enforce that each customer is serviced at least once.  
In \eqref{flow_term} we enforce that for each $l\in \Omega_R$ that the number of LA-arcs leaving a given node (excluding the source,sink) is the same as the number of arcs entering that node.  Since the number of variables in \eqref{my_rmp} is not massive we can solve $\Psi(\cup_{l \in \Omega_R}\Omega_{l,I^-,P^-})$ via solving \eqref{my_rmp}  using any off the shelf LP solver.   Observe that the equality constraint matrix for \eqref{my_rmp} is primal-block angular (grouping by $l$).  Thus we can solve \eqref{my_rmp} using an LP solver that exploits primal block angular structure such as \citep{castro2007interior}.

Observe that solving pricing for the lowest reduced cost $l\in \Omega$ or $l \in\Omega_{\hat{l}}$ does not include dual variables over \eqref{flow_term} as is discussed in \citep{yarkony2021graph,yarkony2022principled,mandal2022local_2}; as these terms would cancel out in the reduced cost of any route.  It is proven in \citep{yarkony2021graph} that since each path from source to sink in a given multi-graph describes a feasible route that $\Psi(\Omega_R,I,P)=\Psi(\cup_{l \in \Omega_R}\Omega_l)$.  .

In Alg \ref{comp_alg} we describe the solution to $\Psi(\Omega)$ with exposition provided below.
\begin{algorithm}[!b]
 \caption{Solving $\Psi(\Omega)$}
\begin{algorithmic}[1] 
\State $\Omega_R\leftarrow$ any initial column (s) or from user \label{init1}  
\State $I^{-},P^- \leftarrow I^{-0},P^{-0}$ \label{init2}  \quad \mbox{ add minimal nodes/LA-arcs needed to express all single customer routes} \quad  
\Repeat{  Loop Solves:  $\Psi(\Omega)$}\label{loop_2} 
\Repeat{  Loop Solves:  $\Psi(\cup_{l \in \Omega_R}\Omega_l)$} \label{outer_2}
\State $x,\pi \leftarrow $ Solve $\Psi(\Omega_R,I^-,P^- )$ via \eqref{my_rmp} \label{step_RMP}
\State $\hat{l}^l \leftarrow \min_{\hat{l} \in \Omega_l}\bar{c}_{\hat{l}}  \quad \forall l \in \Omega_R $ \label{call_ez_pricing}
\State $I^-\leftarrow I^- \cup \cup_{l \in \Omega_R} I_{\hat{l}^l}$ \label{I_aug}
\State $P^-\leftarrow P^- \cup \cup_{l \in \Omega_R} P_{\hat{l}^l}$ \label{P_aug}
\Until{$\bar{c}_{\hat{l}^l}\geq 0  \quad \forall l \in \Omega_R$} \label{do_term_inner}
\State $I^-\leftarrow \{ i \in I^-; \mbox{ s.t. } \exists (l\in \Omega_R,j\in I^-,p \in P^-) \mbox{ s.t. } x^l_{ijp}>0\}$ \quad (OPTIONAL) \label{keep_I}
\State $P^-\leftarrow \{ p \in P^-; \mbox{ s.t. } \exists (l\in \Omega_R,i \in I^-,j\in I^-) \mbox{ s.t. } x^l_{ijp}>0\} $ \quad (OPTIONAL)\label{keep_P}
\State $l^* \leftarrow \min_{l \in \Omega}\bar{c}_l$ \label{call_fancy_pricing}
\State $\Omega_{R} \leftarrow \Omega_{R} \cup l^*$ \label{augment_FR}
\Until{$\bar{c}_{l}^*\geq 0$} \label{loop_2_e}
\end{algorithmic}
\label{comp_alg}
\end{algorithm} 

\begin{enumerate}
    \item Line \ref{init1}:  We receive one or more columns (providing $\beta^l$ terms) from the user. In our implementation initialize $\Omega_R$ with an empty route so that we get an entirely random $\beta^l$ term.
    \item Line \ref{init2}: We initialize $I^-,P^-$ to $I^0,P^0$, which are the minimal nodes/LA-arcs permitting all single customer routes to be created; and are defined below.  
    \begin{subequations}
        \begin{align}
        I^0=\{ (-1,d_0,t_0),(-2,0,0),(u,d_0,\min(t_0-t_{-1u},t^+_u)) \quad   \forall u \in N\}\\
        P^0= \{\cup_{u \in N}(-1,u,\{ \}),(u,-2,\{ \})\}
        \end{align}
    \end{subequations}
        \item Line \ref{loop_2}-\ref{loop_2_e}:  Solve $\Psi(\Omega)$
        \begin{itemize}
        \item Line \ref{outer_2}-\ref{do_term_inner}:  $x,\pi \leftarrow $Solve  $\Psi(\cup_{l \in \Omega_R}\Omega_l)$ 
        \item 
    \begin{itemize}
        \item Line \ref{step_RMP}:  $\pi,x\leftarrow $ Solve the RMP in  \eqref{my_rmp} using an off the shelf solver or by exploiting the primal block angular structure using \citep{castro2007interior}.
        \item 
        Line \ref{call_ez_pricing}:  Compute the lowest reduced cost route in $\Omega_l$ denoted $\hat{l}^l$ for each $l\in \Omega_R$. 
        \item Line \ref{I_aug}-Line \ref{P_aug}:  Add in the nodes/LA-arcs used in any route $\hat{l}^l$ over $l \in \Omega_R$ to the set under consideration.  
        \item Line \ref{do_term_inner}:  Terminate when no route generated has negative reduced cost.  Thus we have solved $\Psi(\cup_{l \in \Omega_R}\Omega_l)$.
        \end{itemize}
        \item Line \ref{keep_I}-\ref{keep_P}:  We remove any nodes/LA-arcs not used in any multi-graph in the solution.   This is optional and can accelerate the solution to the LP in line \eqref{step_RMP} at the expense of increasing the number of inner loop iterations.  We did not do this in our experiments.  
        \item Line \ref{call_fancy_pricing}:  Solve the ERCSPP to generate  $l^*$, which is the lowest reduced cost route in $\Omega$.  This operation is often far more computationally time intensive than a call to price over a family as in Line \ref{call_ez_pricing}.
        
        We should note that we can produce any negative reduced cost column (as long as one exists) without preventing the convergence of optimization.  This is useful early in CG optimization when dual variables may be far from an optimal solution for $\Psi(\Omega)$ \citep{desrosiers2005primer}.
        \item Line \ref{augment_FR}:  Add $l^*$ to $\Omega_R$. 
        \item Line \ref{loop_2_e}:  Terminate optimization when no route has negative reduced cost.  At this point we have solved $\Psi(\Omega)$.  This can be incorporated into a branch-cut-price solver \citep{ropke2009branch}.  We can exploit LA-SRI and LA-rounded capacity inequalities (RCI) as is done in \citep{mandal2022local_2} to tighten the relaxed solution with neglibible additional time taken.  
    \end{itemize}
\end{enumerate}
\section{Pricing Exploiting LA-Routes}
\label{sec_pricing}
In this section we consider the efficient solution to pricing over $\Omega$ and $\Omega_l$ as required in Alg \ref{comp_alg}. The differences between the approaches for pricing over $\Omega$ compared to that over $\Omega_l$ are minor. Thus we describe pricing in this section to be over $\Omega$, except where explicitly noted to correspond to $\Omega_l$.   

\textbf{Summary of this Section:}  Our approach is based on two graphs called the super LA pricing graph and the relaxed graph, which are directed graphs with non-negative weights. The shortest path from source to sink on the super LA pricing graph  corresponds to the lowest reduced cost route in $\Omega$.  In contrast the shortest path from source to sink in the relaxed pricing graph provides a lower bound on the reduced cost of the lowest reduced cost route in $\Omega$. The super LA pricing graph is too large to consider explicitly much less solve pricing over, while the relaxed graph can easily be used to solve pricing due to its much smaller number of nodes and edges. The relaxed graph is constructed by merging nodes in the super LA pricing graph and connecting nodes in the relaxed graph by the lowest cost edge between their parents in the super LA pricing graph. Inspired by \citep{boland2017continuous}  we construct an algorithm that iteratively computes the lowest cost path in the relaxed graph; followed by splitting nodes to ensure that lower bound is tightened \footnote{We introduce what we believe is a new idea here; the notion of capacity discretization as to be shown in \eqref{etaDef}}. We terminate when the route generated by pricing over the relaxed graph provides the lowest reduced cost route in $\Omega$ provably.  We initialize the relaxed graph with the coarsest possible aggregation of nodes so that each member of $N^+$ is associated with exactly one node. We can terminate early when we can project the generated route to a feasible negative reduced cost route. 

We organize this section as follows.  In Section \ref{subsec_orig_PricingGraph} we introduce the super LA pricing graph.  In Section \ref{subsec_relaxedGraph} we introduce the relaxed graph.  In Section \ref{subsec_aug} we describe how to split nodes in the relaxed graph given the lowest cost path in the relaxed graph so as to tighten the relaxation.  In Section \ref{subsec_alg} we consider our algorithm for creating a sufficient relaxed graph so as to generate the lowest reduced cost route in $\Omega$.  In Section \ref{subsec_imp_detail} we consider implementation details.
\subsection{Super LA Pricing Graph}
\label{subsec_orig_PricingGraph}
In this section we describe the super LA pricing graph denoted $H$, which has node set $J$, and edge set $Y$.
There is one node in $H$ for each  $i=(u_i \in N,d_0\geq d_i \geq d_{u_i},t^+_{u_i}\geq t_i\geq t_i^+,M_i\subseteq N-u_i)$.  The source node of $H$ is defined as $(1,d_0,t_0,\{ \})$ while the sink is defined as $(-2,0,0,\{ \})$.  Edges in $Y$ describe the intermediate customers visited (and the associated order) between the customers associated with nodes.  A path from source to a node $i$ (except the sink) indicates a vehicle departing $u_i$ at time $t_i$ with $d_i$ units of demand remaining prior to servicing $u_i$, and that it has serviced all customers in $M_i$ exactly once and no other customers.  A path from source to sink describes a complete route.

Our discussion of edge weights relies on the following two terms, which are $\pi_p$ and $\eta$.  Here $\pi_p$ is the sum of the dual variables associated with $N^+_p$ for a given $p \in P$.  We define $\eta$ s.t. $-\eta*d_0$ is a lower bound on the reduced cost of the lowest reduced cost route in $\Omega$.  Specifically $\eta$ is defined by setting $-\eta$ to be the lowest reduced cost per unit demand consumed over LA-arcs. 
We define $\pi_p$ and $\eta$ below. 
\begin{subequations}
    \begin{align}
     \pi_{p}=\sum_{w \in N^+_p}\pi_w \quad \forall p \in P\\
    \eta=-\min_{p \in P}\frac{c_{p,t^+_{u_p},t^-_{v_p}}-\pi_p}{d_p} \label{etaDef}
    \end{align}
\end{subequations}
%
We construct $Y$ such that edge weights are non-negative, and that for any path corresponding to a route in $l \in \Omega$, that the sum of the edge weights on the path equals $\bar{c}_l+\eta d_0$ .  For each $i\in J,j\in J$ there is non-negative cost to traverse from $i$ to $j$ denoted $\hat{c}_{ij}$.  Typically $\hat{c}_{ij}$ this is infinite, in which case, traveling from $i$ to $j$ is infeasible.  The edges in $Y$ are the $i,j$ pairs for which $\hat{c}_{ij}$ is non-infinite.  We describe $\hat{c}_{ij}$ below using helper terms defined subsequently.
%
\begin{subequations}
\label{super_orig_pricing_graph}
\begin{align}
\hat{c}_{ij}=\eta (d_i-d_j)+\min_{p \in P_{ij}}c_{pt_it_j}-\pi_p \label{hatCij} \quad \forall i \in J,j \in J\\
P^+_{ij}=\{ p \in P; u_i=u_p, u_j=v_p,(d_p=d_i-d_j \quad \mbox{ or } (u_j=-2 \mbox{ and } d_p\leq d_i))\}\\
P_{ij}= \{p \in P^+_{ij}: \quad |M_i \cap (N_p^+\cup v_p)|=0,M_j=M_i\cup N^+_p \}
\end{align}
\end{subequations}
Our helper terms are defined as follows.  
\begin{itemize} 
\item $P^+_{ij}$:  This is the set of LA-arcs that are associated with a given $i,j$ if we ignore consistency with regards to the customers already visited.  If $j$ is not the sink then the demand serviced on the LA-arc must equal $d_i-d_j$. However if $j$ is the sink then any amount of demand less than or equal to $d_i$ can be serviced on the LA-arc.  This is done so as to permit the expression of routes servicing less than $d_0$ units of demand.
\item $P_{ij}$:  This is the subset of $P^+_{ij}$ that is consistent with regards to customers already visited. When $P_{ij}$ is empty then $\hat{c}_{ij}=\infty$, meaning that the edge  $(i,j)$ does not exist in $Y$.  
\end{itemize}
Observe that $\hat{c}_{ij}$ is non-negative  by definition of $\eta$  in \eqref{etaDef}.  The minimizing $p$ in \eqref{hatCij} is described as $p^*_{ij}$, and the ordering of the customers visited in $N_p^+$ (where $p=p^*_{ij}$) is denoted $N^+_{ij}$.   
 Traversing the edge from $i$ to $j$ indicates that the LA-arc $p^*_{ij}$ is used in the solution to pricing.  A large subset of routes in $\Omega$ can be represented by paths from source to sink in $H$; and all paths from source to sink are associated with feasible routes in $\Omega$.  We should note that there are exactly two classes of routes in $\Omega$ can not represented in $H$.  However neither of these classes of routes can contain the lowest reduced cost route in $\Omega$, and thus need not be expressed $H$. The first class contains routes that have a sub-optimal orderings of customers with regards to cost.  The second are routes that would use an LA-arc on edge $i,j$ for which the associated LA-arc $p$ satisfies $\hat{c}_{ijp}>\hat{c}_{ij}$. Such an arc would not be used in an optimal solution to pricing. Observe that the following four properties are satisfied for any route $l$ corresponding to a path from source to sink in $H$.  
%
%
\begin{itemize}
\item $\eta$ Correctness: Observe that no matter how much demand is serviced in route $l$ that $\eta \sum_{ij \in Y^l}(d_i-d_j)=d_0\eta$; where the edges visited in $l$ are denoted $Y^l$.  This means that the total cost of edges on the route is $\bar{c}_l+\eta d_0$.
\item Demand Feasibility:  Observe that $\sum_{u \in N}d_u a_{ul}\leq d_0$.  This means that no more demand is serviced than $d_0$.
\item Temporal Correctness:  Observe that  $t^+_u\geq t^l_u\geq t^-_u$ for all $u\in N$ s.t. $a_{ul}=1$; where $t^l_u$ is the time when $u$ is departed from. Thus all time windows are respected. 
\item Elementarity (No Cycles):  Each customer is included no more than once in the route. Thus the route is elementary.
\end{itemize}
When pricing over a family $\Omega_l$, as required in Line \ref{call_ez_pricing} of Alg \ref{comp_alg},  we only use valid LA-arcs for the family of $l$ when constructing $H$.  Thus we replace $P^+_{ij}$ with $P^{l+}_{ij}$ in the definition of $P_{ij}$ which is defined as follows (as is expressed analogously for the master problem in \eqref{def_arcs_family}).  
\begin{align}
\label{pdefFrom}
P^{l+}_{ij}=\{ p \in P^+_{ij}, \beta^l_{u_p}<\beta^l_{v_p}, \beta^l_{u_p}<\beta^l_{w}<\beta^l_{v_p} \quad \forall w \in N_p\}
\end{align}
We refer to the cost of the lowest cost path from the source to sink in $H$ as $\hat{c}_{opt}$. 
\subsection{Relaxed Graph}
\label{subsec_relaxedGraph}
In this section we describe the relaxed graph denoted $H^-$, which has node set $J^-$, and edge set $Y^-$.  Each node $g \in J^-$ is defined as $g=(u_g,d^-_g,d^+_g,t^-_g,t^+_g,M^-_g,M^+_g)$. Each node $g\in J^-$ is associated with a subset of the nodes in $J$.  The nodes in $J$ are partitioned between $J^-$ node sets.  Thus $\cup_{g \in J^-}J_g=J$ and $|J_{g}\cap J_f|=0$ for all $f\neq g$.  
The source and the sink are respectively written as $\{ (-1,d_0,d_0,t_0,t_0,\{ \},\{ \}),\\(-2,0,d_0,0,t_0,\{ \},\{ N\})$.
We define $J_g \subseteq J$ as follows.  
\begin{subequations}
\begin{align}
J_g=\{ i \in J,  u_g=u_i, \\
d^-_g\leq d_i \leq d_g^+, \\
t^-_g\leq t_i \leq t_g^+, \\
M^-_g\subseteq M_i \subseteq M^+_g \} 
\end{align}
\end{subequations}
For each $f\in J^-,g \in J^-$ we connect $f$ to $g$ in $H^-$ with the lowest cost edge between any $i \in J_f, j\in J_g$ ; if such an edge exists and otherwise do not connect $f$ to $g$.
We write the associated edge weight as $\hat{c}_{fg}$, which we express as follows.  
\begin{align}
\label{hatgfDef}
\hat{c}_{fg}=\min_{\substack{i \in J_f\\ j\in J_g}}\hat{c}_{ij} \quad \forall f \in J^-,g \in J^-
\end{align}
We write the efficient computation of $\hat{c}_{fg}$ below, using helper terms defined inside the equation set with further exposition below the equation set.
\begin{subequations}
\label{edge_vec}
\begin{align}
\hat{c}_{fg}=\min_{p \in P_{fg}}\hat{c}_{fgp} \label{edge_vec0}\\
\hat{c}_{fgp}=\eta d_{fgp}-\pi_{p}+c_{p,t_f^+,t_g^-}\\
d_{fgp}=\max(d_f^- -d^+_g[u_g\neq -2],d_p)\\
P^+_{fg}=\{ p \in P; u_f=u_p, u_g=v_p,d_f^+-d_p\geq d_g^-;
d_f^- -d_p\leq d^+_g\}\\
P_{fg}= \{p \in P^+_{fg}: \quad |M^-_f \cap (N_p^+ \cap v_p)|=0,M^-_g \subseteq( M_f^+\cup N_p^+),  \quad (M_f^-\cup N^+_p) \subseteq M^+_g \}
\end{align}
\end{subequations}
\begin{itemize}
\item $\hat{c}_{fgp}:  $ This is the reduced cost of the lowest reduced cost path associated with LA-arc $p$, leaving $u_f$ at $t^+_f$ and leaving $u_g$ at $t^-_g$; plus it is augmented by the $\eta$ times amount of demand serviced.  
Observe that $c_{p,t_1,t_2}$ is non-increasing as $t_1$ increases and non-increasing as $t_2$ decreases.  This is because increasing $t_1$ and or decreasing $t_2$ provides only more variability to the set of possible orderings of customers to use for LA-arc $p$.  Hence we can use $c_{p,t_f^+,t_g^-}$ instead of trying all possible values of $t^+_f\geq t_1\geq t_f^-$,$t_g^+\geq t_2\geq t^-_g$.
\item $d_{fgp}:  $  This is the amount of demand serviced on arc $f,g$ (or minimum demand remaining at $f$ if this is the last LA-arc in the route and $d_f^->d_p$).
\item $P^+_{fg}$:  This is the set of LA-arcs consistent with regards to demand remaining.  When pricing over a family $\Omega_l$ instead of $\Omega$ we  replace $P^+_{fg}$ with $P^{l+}_{fg}$ which is defined below.  
\begin{align}
\label{pdefFrom2}
P^{l+}_{fg}=\{ p \in P^+_{fg}, \beta^l_{u_f}<\beta^l_{v_g}, \beta^l_{u_f}<\beta^l_{w}<\beta^l_{u_g} \quad \forall w \in N_p\}
\end{align}
 Observe that \eqref{pdefFrom2} is analogous to \eqref{pdefFrom} where we simply replace $i,j$ with $f,g$ respectively.
\item $P_{fg}$ is the subset of $P^+_{fg}$ consistent with the customers visited thus far.   
\end{itemize}
For a given $f,g$ the minimizing $p$ in \eqref{edge_vec0} is denoted $p^*_{fg}$ and the associated ordering of customers in $N^+_{p^*_{fg}}$ is denoted $N^+_{fg}$.
The shortest path in $H^-$ from source to sink is easily computed via Dykstra's algorithm, since all weights in $Y^-$ are non-negative.  We describe this shortest path as $l^{*-}$ and use $l\leftarrow l^{*-}$ to denote the ``column" associated with this path; where $l$ need not lie in $\Omega$, and it need not obey the four properties of $\eta$ Correctness, Demand Feasibility, Temporal Correctness, or Elementarity. 
 However observe that if  $l^{*-}$ did obey these four properties then $l$ would describe the lowest reduced cost column in $\Omega$. Let $Y^{-l}$ denote the edges in $Y^-$ used in path $l$. 
The cost of the path for $l$, which is denoted $\hat{c}_{Lb}$ is defined to be the sum of the edges on the path; thus $\hat{c}_{Lb}=\sum_{fg \in Y^{-l}}\hat{c}_{fg}$. 
Observe that since $H^-$ assigns costs to edges, which are lower bounds to the corresponding costs assigned to edges in $H$, then $\hat{c}_{opt}\geq \hat{c}_{Lb}$.   
 We define the ordering of nodes on the relaxed graph taken by $l$, from first to last, as $[g^l_0...g^l_{k},g^l_{k+1}...g^l_{|g^l|}]$ where $|g^l|$ is the number of such nodes.  We use $[-1,u^l_{1},u^l_2....-2]$ to denote the customers/depots visited in route $l$ ordered as described by the $N^+_{fg}$ terms.  For any $k$ we use $q^l_k$ to denote the number of customers serviced prior to reaching $g^l_{k}$ in route $l$.  For example for a route $[-1,8,3,5,6,7,9,-2]$ with nodes corresponding the positions of customers/depots $-1,8,7,-2$ then $q^l_0=q^l_1=0$,$q^l_2=4$ and $q^l_3=6$.  
\subsection{Augmenting the Relaxed Graph by Splitting Nodes}
\label{subsec_aug}
In this section we consider the classes of mechanisms causing $\hat{c}_{Lb}$ to be a loose bound on $\hat{c}_{opt}$ meaning $\hat{c}_{opt}>\hat{c}_{Lb}$.  We first describe these mechanism classes.  Next we present a sufficient condition to ensure that a particular mechanism can not occur for a given route $l$. Then we demonstrate how to alter the relaxed graph $H^-$ so as to prevent a particular mechanism from being active by splitting nodes in $J^-$.  The source and the sink are never split.  
\begin{itemize}
\item \textbf{Case One:} $\eta$ Correctness is Violated:
\begin{itemize}
\item  \textbf{Description:} In this case $\eta$ is not added to the cost $d_0$ times.  By this we mean that $\sum_{fg \in Y^{-l}} d_{fgp^*_{fg}} < d_0$.  
\item \textbf{New Notation Required:}  We use $d_{0:k}$ be the total capacity remaining immediately prior to reaching node $g^l_k$ meaning $d_{0:k}=d_0-\sum_{q^l_k\geq i\geq 1}d_{u^l_{i}}$.
\item 
\textbf{Sufficient Condition for Avoidance:  }This is avoided when the following is satisfied for each $k$ s.t. $|g^l|>k\geq 1$:  $d^-_{g^l_k}=d_{0:k}$.   
\item 
\textbf{Correction by Node Splitting}: For each $k$ for which $|g^l|>k\geq 1$  and  $d^+_{g^l_k}\geq d_{0:k}>d^-_{g^l_k}$ we replace node $g^l_k$ in $J^-$ with two nodes $g_1,g_2$ that are identical to $g^l_k$ with the following exceptions. 
\begin{subequations}
\label{rule1}
\begin{align}
d^+_{g_1}=d_{0:k}-1\\
d^-_{g_2}=d_{0:k}
\end{align}
\end{subequations}
\end{itemize}
\item \textbf{Case Two:}   Capacity Exceeded
\begin{itemize}
\item \textbf{Description:}  In this case the total demand serviced in the route exceeds the capacity $d_0$ meaning that $\sum_{fg \in Y^{-l}}d_{p^*_{fg}}>d_0$.  
\item \textbf{Sufficient Condition for Avoidance:  } 
This is avoided when the following is satisfied for each $|g^l|>k\geq 1$ :  $d^+_{g^l_k}=d_{0:k}$
\item \textbf{Correction by Node Splitting}: For each $k$ for which $|g^l|>k\geq 1$  and  $d^+_{g^l_k}> d_{0:k}\geq d^-_{g^l_k}$ we replace node $g^l_k$ in $J^-$ with two nodes $g_1,g_2$ that are identical to $g^l_k$ with the following exceptions.
\begin{subequations}
\label{rule2}
\begin{align}
d^-_{g_1}=d_{0:k}+1\\
d^+_{g_2}=d_{0:k}
\end{align}
\end{subequations}
\end{itemize}
\item \textbf{Case Three:}  Time Window Infeasibility
\begin{itemize}
\item \textbf{Description:} In this case we are not able to visit all of the customers in the order associated with route $l$.
\item \textbf{New Notation Required}:   We use $\hat{t}_{0:k}$ to denote the correct departure time (with no unnecessary waiting) at the customer associated with $g^l_k$.  We define this recursively below using helper term $t_{0:i}$ to denote the correct departure time (with no unnecessary waiting) at the customer associated with $u^l_i$.
\begin{subequations}
\begin{align}
\hat{t}_{0:k}=t_{0:q^l_k+1}\\
t_{0:0}=t_0\\
t_{0:i}=\min(t^+_{u^l_i},t_{0:i-1}-t_{u^l_{i-1}u^l_{i}}])-\infty[t_{0:i-1}-t_{u^l_{i-1}u^l_{i}}<t^-_{u^l_i}]
\end{align}
\end{subequations}
\item \textbf{Sufficient Condition for Avoidance:  }
This is avoided when the following is satisfied for each $k$ s.t $|g^l|>k\geq 1$:   $t^+_{u_g}=\hat{t}_{0:k}$, where where $g=g^l_k$. 
\item \textbf{Correction by Node Splitting}:   For each $k$ for which $|g^l|>k\geq 1$ and $t^+_{g^l_k}> \hat{t}_{0:k}\geq t^-_{g^l_k}$ we replace node $g^l_k$ in $J^-$ with two nodes $g_1,g_2$ that are identical to $g^l_k$ with the following exceptions.

\begin{subequations}
\label{rule3}
\begin{align}
t^-_{g_1}=\hat{t}_{0:k}+1\\
t^+_{g_2}=\hat{t}_{0:k}
\end{align}
\end{subequations}
\end{itemize}
\item \textbf{Case Four:}  Repeated Customers
\begin{itemize}
\item \textbf{Description:}   There exists a cycle of customers in the route.  Note that when pricing over some family $\Omega_l$ this can not occur and hence this possibility need not be considered.
\item \textbf{New Notation}:  We define the minimum length cycle in $l$ as the cycle of customers for which the number of intermediate customers associated with nodes $g^l_k$ is minimized.  This is described as follows.  Given that the cycle starts and ends and positions $i_1,i_2$ with customer $u$ where $u=u^l_{i_1}=u^l_{i_2}$ let $k_1,k_2$ be the nodes that most tightly enclose the cycle.  Thus $k_1,k_2$ are the largest/smallest values respectively for which  $q^l_{k_1}+1\leq i_1  <q^l_{k_{2}-1}+1<i_2\leq q^l_{k_2}+1$.  The intermediate nodes are the nodes $k_3$ for which $k_2>k_3>k_1$.  When multiple minimum length cycles exist one is selected arbitrarily.  
\item \textbf{Sufficient Condition for Avoidance:  }
This cycle can only occur if there exists a $k_3$ for which $k_2>k_3>k_1$ and $g=g^l_{k_3}$ satisfies $u\in M^+_{g}$ and $u \notin M^-_{g}$.  
\item \textbf{Correction by Node Splitting}: For each $k_3$ s.t $k_2>k_3>k_1$ for which, $u \in M^+_{g^l_{k_3}}$ and $u\notin M^-_{g^l_{k_3}}$ we replace node $g^l_{k_3}$ (denoted $g$) in $J^-$ with two nodes $g_1,g_2$ that are identical to $g$ with the following exceptions.
\begin{subequations}
\label{rule4}
\begin{align}
M^-_{g_1}=M^-_{g} \cup u\\
M^+_{g_2}=M^+_{g} - u
\end{align}
\end{subequations} 

\end{itemize}
\end{itemize}%

\subsection{A Node Splitting Algorithm for Efficient Pricing }
\label{subsec_alg}

In this section we describe our algorithm for solving for the lowest reduced cost route by constructing sufficient set $J^-$ s.t the lowest cost path in $H^-$ violates none of the cases from Section \ref{subsec_aug}.  This algorithm iterates between \textbf{(1)}solving for the lowest cost path in $H^-$ and \textbf{(2)}identifying at least one violation and correcting it using node splitting as described in Section \ref{subsec_aug}.  This algorithm terminates when no such violations exist.  

We found it effective to correct one violation at each iteration.  Specifically we prioritized operations on cases (1),(2),(3),(4) in that order. Alternative orderings can however be used to still produce an efficient pricing procedure. We see that case (1) is intuitively most important since not having $\eta d_0$ be a component in $\hat{c}_{Lb}$ can make $\hat{c}_{Lb}$ much less than $\hat{c}_{opt}$.  Similarly we observe that cycles of customers often disappear once feasibility with regards to demand/time is enforced so we put case (4) last. We use Viol($l^{*-})$ to denote the violation selected.  
\begin{algorithm}
\caption{Solving Pricing Efficiently}
\begin{algorithmic}[1]
\State $J^- \leftarrow  J^{0}$
\label{init_I_minus}
\Repeat \label{loopStartPricing}
\State Compute $\hat{c}_{f,g}$ for each $f \in J^-,g \in J^-$. \label{update_weights}
\State $l^{*-} \leftarrow $ Shortest Path in $H^-$ \label{get_short_path}
\Switch{Viol($l^{*-}$)} \label{startSwitch}
    \Case{1:  }
      $J^-\leftarrow $ Split using \eqref{rule1} \label{split1}
    \EndCase
    \Case{2:  }
      $J^-\leftarrow $ Split using \eqref{rule2} \label{split2}
    \EndCase
    \Case{3:  }
      $J^-\leftarrow $ Split using \eqref{rule3} \label{split3}
    \EndCase
    \Case{4:  }
      $J^-\leftarrow $ Split using \eqref{rule4}  \label{split4}
    \EndCase
  \EndSwitch \label{endSwitch}
\Until{Viol($l^{*-}$) is empty} \label{TermLoopPrice}
\State Return $l^{*-}$ \label{returnPriceSol}
\end{algorithmic}
\label{fastPricing}
\end{algorithm}

\begin{itemize}
    \item Line  \ref{init_I_minus}:  Initialize $J^-$ to be the coarsest partition of nodes denoted $J^{0}$, which we define as follows.
\begin{subequations}
\begin{align}
    J^{0}=\{ (-1,d_0,d_0,t_0,t_0,\{ \},\{ \}),\\(-2,0,d_0,0,t_0,\{ \},\{ N\}),\\ 
    (u,d_u,d_0,t^-_u,t^+_u,\{ \},N-u) \forall u \in N \}
\end{align}
\end{subequations}
    \item Line  \ref{loopStartPricing}-\ref{TermLoopPrice}:  Iterate between \textbf{(1)} generating the lowest cost path from source to sink in $H^-$ and \textbf{(2)} augmenting $J^-$ until that path describes a column that lies in $\Omega$,  and the total cost of the incurred $\eta$ terms is $d_0\eta$.  When this is satisfied then the path corresponds to a column that is the lowest reduced cost column in $\Omega$.
    \begin{itemize}
    \item Line \ref{update_weights}:  Compute edge weights.  This need only be done for all edges not present in the previous iteration of this loop. 
    \item Line \ref{get_short_path}:  We compute the lowest cost path from source to sink, which we denote as $l^{*-}$. 
    \item Line \ref{startSwitch}-\ref{endSwitch}:  If a violation has exists in $l^{*-}$ then we split some nodes in $J^-$ so as to correct the violation.  
    \begin{itemize}
    \item Line \ref{split1}:  Split nodes so as to enforce that the total amount of $\eta$ terms used in $l^{*-}$ is $d_0 \eta$.
    \item Line \ref{split2}:  Split nodes so as to enforce that the total demand serviced in $l^{*-}$ does not exceed $d_0$.
    \item Line \ref{split3}:  Split nodes so as to enforce time window feasibility in $l^{*-}$.  
    \item Line \ref{split4}:  Split nodes so as to remove a cycle of customers found in $l^{*-}$ from feasibility.  
    \end{itemize}
    \item Line \ref{TermLoopPrice}:  Terminate when $l^{*-}$ lies in $\Omega$ and the total weight of $\eta$ terms included is $d_0\eta$
    \end{itemize}
    \item Line:  \ref{returnPriceSol}:  Return $l^{*-}$, which is the lowest reduced cost route in $\Omega$. 
\end{itemize}
\subsection{Implementation Details}
\label{subsec_imp_detail}
We now describe the following implementation details for Alg \ref{fastPricing}.  
\begin{itemize}
\item
\textbf{Retaining Old Nodes for Exact Pricing}:  We seek to improve the efficiency of Alg \ref{fastPricing} over $\Omega$ by recycling $J^-$ from previous iterations. 
Thus we initialize  $J^-$ with the $J^-$ produced from the previous call to pricing over $\Omega$ using Alg \ref{fastPricing}.  We only initialize $J^-$ with $J^0$ on the first call to pricing over $\Omega$.  We could remove nodes from $J^-$ that are not used, in an analogous way to Lines \ref{keep_I}-\ref{keep_P} of Alg \ref{comp_alg} though did not explore this experimentally.  
\item
\textbf{Retaining Old Nodes for Pricing over Families}:   We seek to improve the efficiency of pricing over families by recycling $J^-$ terms.  We refer to the nodes in $J^-$ in the last call to pricing over $\Omega_l$ as $J^{-l}$.  We initialize $J^-$ as $J^{-l}$ when pricing over $\Omega_l$.  The set $J^{-l}$ is thus augmented after each call to pricing over the family $\Omega_l$.  We initialize $J^{-l}$ as $J^{0}$ prior to the first call to pricing over $\Omega_l$.  As when pricing over $\Omega$ we could remove nodes from $J^{-l}$ that are not used in an analogous way to Lines \ref{keep_I}-\ref{keep_P} of Alg \ref{comp_alg} though did not explore this experimentally.
\item
\textbf{Early Termination of Exact Pricing:  }  It is common in the CG literature to stop pricing early when a negative reduced cost column has been generated.  This is based on the observation of the ``heading in effect"\citep{desrosiers2005primer}, which states that early in the course of CG that the dual variables are far from their final values, and hence convergence of CG is not vastly accelerated by using the lowest reduced cost column instead of any negative reduced cost column.  Thus we may choose to terminate pricing when we can project the solution $l^{*-}$ generated on \eqref{get_short_path} of Alg \ref{fastPricing} to a column $l^*\in \Omega$ with negative reduced cost within a tolerance of optimality.  Thus for a for a factor of tolerance $\alpha \in (0,1)$ we return $l^*$ if $\bar{c}_{l^*}<\alpha (\hat{c}_{Lb}-\eta d_0)$.  We use $\alpha=0.2$ in our experiments when pricing over $\Omega$.  

We now consider the projection of $l^{*-}$ to $l^*$.  We initialize $l^*$ as containing $[-1,-2]$. We iterate over the customers of $l^{*-}$ in order of when visited in $l^{*-}$ from first to last. Given a customer $u$, we add $u$ just before $-2$ in $l^*$ if $l^*$ remains feasible given this addition. Here feasible indicates that $l^*$ lies in $\Omega$. 
\end{itemize}
\section{Local Area Arcs with Time Windows}
\label{sec_la_time}
In this section we describe the efficient computation of $c_{p,t_1,t_2}$ as required in the previous sections.  This is done by exploiting the structure of a dynamic program so that we need not consider all orderings of customers in all LA-arcs.  This dynamic program produces a structure so that given any $p,t_1,t_2$ a very fast computation can be done to determine $c_{p,t_1,t_2}$.  This exploits the fact that for a given $p$ there are only a small number of orderings of the intermediate customers that are optimal for any pair of times.  This also exploits the joint computation of all terms associated with LA-arcs so as to avoid wasting computation.  

We organize this section as follows.  In Section \ref{subsec_make_term_notation} we provide notation needed to express the material in this section.  In Section \ref{subsec_make_term_recur} we describe a recursive relationship for the cost of an LA-arc $p$ given  earliest departure time from $u_p$ and latest departure time at $v_p$.  In Section \ref{subsec_make_term_fronteir} we introduce the concept of an efficient frontier for each $p\in P$ trading off latest departure time from $u_p$ (later is preferred); earliest departure time from $u_p$ when no waiting is required (earlier is preferred); and cost (lower cost is preferred).  In Section \ref{subsec_make_term_alg} we provide a dynamic programming based algorithm to compute the efficient frontier for all LA-arcs jointly.  
\subsection{Notation for LA-Arcs Modeling Time}
\label{subsec_make_term_notation}
In this section we provide additional notation required for discussion of LA-arcs modeling time.  
Let $P^a$ be a super-set of $P$ defined as follows.  Here $p=(u_p,v_p,N_p)$ for $u_p \in N,v_p\in N^+,N_p\subseteq N$ lies in $P^a$ if any only if the following condition holds. There exists a $u\in N$ s.t. $u_p \in (u \cup N_u),v_p\in N^+-(N_u\cup u), N_p \subseteq N_u$.  We use $R_p$ to denote the set of all possible paths starting at $u_p$ ending at $v_p$ and visiting $N_p$ in some order. We describe a path $r$ as an ordered sequence of customers listed from beginning to end as follows $[u^r_1,u^r_2,u^r_3,...v^r]$. We also define a path $r^-$, to be the same as $r$ except with $u^r_1$ removed meaning $r^-=[u^r_2,u^r_3,...v^r]$.

For any given time $t$ and path $r\in R_p$, we define $T_{r}(t)$ as the earliest time a vehicle could depart $v^r$, if that vehicle departs $u^r_1$ with $t$ time remaining and follows the path of $r$.  We write $T_r(t)$ mathematically below using $-\infty$ to indicate infeasibility. 
\begin{align}
\label{trtDef}
    T_{r}(t)=T_{r^-}(-t_{uw}+\min(t,t^+_u))-\infty*[t<t^-_u] \quad \forall r\in R_p,p \in P, u^r_1=u,w=u^r_2
\end{align}
Given $T_{r}(t)$ we express $c_{p,t_1,t_2}$ as follows.
\begin{align}
\label{lookupTravelTime_orig0xx}
c_{p t_1 t_2}=\min_{\substack{r \in R_p\\ T_{r}(t_1)\geq t_2}}c_{r} \quad \forall p\in P,t_1,t_2
\end{align}
Here $c_{p t_1 t_2}=\infty $ if no $r$ exists satisfying $T_{r}(t_1)\geq t_2$.  Observe that  using \eqref{lookupTravelTime_orig0xx} to evaluate $c_{p t_1 t_2}$ is challenging as we would have to repeatedly evaluate $T_{r}(t)$ in a nested manner.  
\subsection{Recursive Definition of Departure Time}
\label{subsec_make_term_recur}
In this section we provide an alternative characterization of $T_r(t)$, that later provides us with a mechanism to efficiently evaluate $c_{p t_1 t_2}$.  For any $r\in R_p,p\in P^a$ we describe $T_{r}(t)$ using the helper terms $\tau_{r1},\tau_{r2}$ defined as follows.   %
We define $\tau_{r1}$ to be the earliest time that a vehicle could leave $u^r_{1}$ without waiting at any customer in $N^+_p\cup v_p$ if $t^-$ terms were ignored (meaning all $t^-$ terms are set to $-\infty$). Similarly we define $\tau_{r2}$ as the latest time a vehicle could leave $u^r_1$ if $t^+$ terms were ignored (meaning all $t^+$ terms are set to $\infty$).   
Below we define $\tau_{r1},\tau_{r2}$ recursively.
\begin{subequations}
\label{tauDef}
\begin{align}
\tau_{r1}=\min(t^+_u,t^+_v+t_{uv}) \quad \forall r \in R_p, p=(u,v,\{ \}),p \in P^a\\
    \tau_{r2}=\max(t^-_u,t^-_v+t_{uv}) \quad \forall r \in R_p, p=(u,v,\{ \}),p \in P^a\\
    \tau_{r1}=\min(t^+_u,\tau_{r^-1}+t_{uw}) \quad \forall r \in R_p, u^r_{1}=u,u^r_{2}=w, |N_p|>0, p \in P^a\\ 
    \tau_{r2}=\max(t^-_u,\tau_{r^-2}+t_{uw}) \quad \forall r \in R_p,u^r_{1}=u,u^r_{2}=w, |N_p|>0, p \in P^a 
\end{align}
\end{subequations}
Let $c_r$ denote the total travel distance on path $r$.  We rewrite $T_{r}(t)$ using $\tau_{r1},\tau_{r2}$, with intuition provided subsequently (with proof of equivalence in Appendix \ref{sec_proof}).   
\begin{align}
\label{ituDef}
    T_{r}(t)=-c_r+\min(t,\tau_{r1}) -\infty[\min(t,t^+_{u_p})<\tau_{r2}] \quad \forall p \in P^a,r \in R_p 
\end{align}
We now provide an intuitive explanation of \eqref{ituDef}.  Observe that leaving $u^r_1$ after $\tau_{r2}$ indicates that the path is infeasible. 
 Observe that given that $r$ is fixed that leaving $u^r_1$ at time $t$ for $t> \tau_{r1}$ incurs a waiting time of $t-\tau_{r1}$ over the course of the path and otherwise incurs a waiting time of zero.  We then subtract off the travel time $c_r$ to gain the departure time at $v^r$ if the path is feasible.  Given $R_p$ and \eqref{ituDef}  we can write $c_{p,t_1,t_2}$ as follows.  
\begin{align}
\label{lookupTravelTime_orig0}
c_{p t_1 t_2}=\min_{\substack{r \in R_p\\ t_1\geq \tau_{r2}\\ t_2\leq -c_r+\min(t_1,\tau_{r1})}}c_{r} \quad \forall p\in P,t_1,t_2
\end{align}
If no $r$ satisfies the constraints in \eqref{lookupTravelTime_orig0} then $c_{p t_1 t_2}=\infty$.  Observe that $|R_p|$ grows factorially with $|N_p|$ thus using \eqref{lookupTravelTime_orig0} is impractical when $|N_p|$ grows.
\subsection{Efficient Frontier}
\label{subsec_make_term_fronteir}
In this section we describe a sufficient subset of $R_p$ denoted $R^*_p$ s.t. that applying \eqref{lookupTravelTime_orig0} over that subset produces the same result as if all of $R_p$ were considered.  
Given $p$ a necessary criteria for $r \in R_p$ being the lowest cost path for some times $t_1,t_2$ is that $r$ is not Pareto dominated by any other such path in $R_p$ with regards to latest time to start the arc (later is preferred), cost (lower is preferred), and earliest time to start the arc without waiting (earlier is preferred).  Thus we prefer smaller values of $\tau_{r2},c_r$ 
and larger values of $\tau_{r1}$.  
We use $R^*_p\subseteq R_p$ to denote the efficient frontier of $R_p$ with regards to $\tau_{r2},c_r,\tau_{r1}$. Here $r \in R^*_p$ if no $\hat{r}\in R_{p}$ exists for which the following inequalities hold with (and at least one is strict).  
    \textbf{(1)}$\tau_{r2}\geq \tau_{\hat{r}2}$.  
    \textbf{(2)}$c_{r}\geq c_{\hat{r}}$.    
    \textbf{(3)}$\tau_{r1}\leq \tau_{\hat{r}1}$.  We use the following stricter criteria for \textbf{(3)}: there is no time $t$ when we leave $u_p$ using $r$ in which we could depart at $v_p$ before we could depart $v_p$ using $\hat{r}$. This criteria is written as follows $\tau_{r1}-c_{r}\leq \tau_{\hat{r}1}-c_{\hat{r}}$.  Observe that given the efficient frontier $R^*_p$ that we can compute $c_{pt_1t_2}$ as follows.  
\begin{align}
\label{lookupTravelTime}
c_{p t_1 t_2}=\min_{\substack{r \in R^*_p\\ t_1\geq \tau_{r2}\\ t_2\leq -c_r+\min(t_1,\tau_{r1})}}c_{r} \quad \forall p\in P,t_1,t_2
\end{align}
If no $r$ satisfies the constraints in \eqref{lookupTravelTime} then $c_{p t_1 t_2}=\infty$.  We sort $R^*_p$ by $c_r$ from smallest to largest so that we can evaluate \eqref{lookupTravelTime} without necessarily testing each member of $R^*_p$.  Observe that via \eqref{lookupTravelTime} that for any subset of $R^*_p$ with identical $\tau_{r1},\tau_{r2},c_{r}$ terms we need only retain one such term since they would produce identical results.  Observe that we can not construct $R^*_p$ by removing elements from $R_p$ as enumerating $R_p$ would be too time intensive. 
\subsection{An Algorithm to Generate the Efficient Frontier}
\label{subsec_make_term_alg}

In this section we provide an efficient algorithm to construct all $R^*_p$ jointly without explicitly enumerating $R_p$.  
In order to achieve this we exploit the following observation.  Observe that if $r\in R^*_p$ then $r^-$ must lie in $R_{\hat{p}}^*$ where $\hat{p}=(u_{\hat{p}}=u^r_2,v_{\hat{p}}=v^r,N_{\hat{p}}=N_{p}-u^r_2) $. 

We construct $R^*_p$ terms by iterating over $p\in P^a$ from smallest to largest  with regard to $|N_p|$, and constructing $R^*_p$ using the efficient frontier from the previously generated efficient frontiers.  For a given $p$ we first add $u_p$ in front of all possible $r^-$; to construct a candidate list; then remove all $r$ from the candidate list that do not lie on the efficient frontier.  Observe that the base case where $|N_p|$ is empty then there is only one member of $R_p$ which is defined by sequence $r=[u_p,v_p]$ ; when feasible and otherwise is empty. In Alg \ref{better_ver} we describe the construction of $R^*_p$, which we annotate below.
\begin{algorithm}
\caption{Computing the Efficient Frontier for all $p\in P^a$ }
\begin{algorithmic}[1]
\For{$p\in P^a$ from smallest to largest in terms of $|N_p|$ with $|N_p|>0$} \label{outer_0}
\State $R^*_{p}\leftarrow \{ \}$ \label{init_rp}
\For{$w \in N_p; \hat{p}=(w,v_p,N_p-w),r^- \in R^*_{\hat{p}}$} \label{outer1}
\State $r\leftarrow [u_p,r^-]$ \label{grab_gplus}
\State Compute $c_{r}$, and $\tau_{r1}$, $\tau_{r2}$ via \eqref{tauDef} \label{comp_tau}
\If{$\tau_{r2} \le t^+_u$} \label{if_add}
\State $R^*_{p}\leftarrow R^*_p \cup r$  \label{augment_rp}
\EndIf \label{if_add_e}
\EndFor \label{outer1_e}
\For{$r \in R^*_p,\hat{r} \in R^*_p-r$}\label{rem_dominated}
\If{$c_{r}\geq c_{\hat{r}}$ and $\tau_{r2}\geq \tau_{\hat{r}2}$ and $\tau_{r1}-c_{r}\leq \tau_{\hat{r}1}-c_{\hat{r}}$ } \label{if_0}
\State $R^*_{p}\leftarrow R^*_{p}-r $ \label{subtrTerm} 
\EndIf \label{if_e}
\EndFor \label{rem_dominated_e}
\EndFor \label{outer_0_e}
\end{algorithmic}
\label{better_ver}
\end{algorithm}
\begin{itemize}
    \item Line \ref{outer_0}-\ref{outer_0_e}:  Iterate over $p \in P^a$ from smallest to largest with regards to $|N_p|$ and given $p$ compute $R^*_p$.
    \begin{itemize}
    \item Line \ref{init_rp}: Initialize $R^*_p$ to be empty.
    \item Line \ref{outer1}-\ref{outer1_e}:  Compute all possible terms for $R^*_p$ by adding a customer $u_p$ in front of all potential predecessor $r^-$ terms.
    \begin{itemize}
        \item Line \ref{grab_gplus}:  Add $u_p$ to the front of $r^-$ creating $r$.
        \item Line \ref{comp_tau}:  Compute $c_{r}$, and  $\tau_{r1}$, $\tau_{r2}$ via \eqref{tauDef}.
        \item Line \ref{if_add}-\ref{if_add_e}:  If $r$ is feasible then we add $r$ to $R^*_p$.
    \end{itemize}
    \item Line \ref{rem_dominated}-\ref{rem_dominated_e}:  Iterate over members of $R^*_p$ and remove terms that do not lie in the efficient frontier.   
    \begin{itemize}
    \item Line \ref{if_0}-\ref{if_e}:  If $r$ is dominated by $\hat{r}$ then we remove $r$ from $R^*_p$.  
    \end{itemize}
\end{itemize}
\end{itemize}
We compute $R^*_p$ for all $p$ using Alg \ref{better_ver} once prior to running Alg \ref{comp_alg}.  We use \eqref{lookupTravelTime} on demand to compute $c_{pt_1,t_2}$ as required  over the course of Alg \ref{comp_alg}, and Alg \ref{fastPricing}. 
\section{Experimental Validation}
\label{sec_exper}
In this section we demonstrate the value of GM and LA-arcs independently and jointly.  To this end we compare the performance of our approach as we alter the number of LA-neighbors and apply GM vs standard CG.  Here standard CG indicates that we use Alg \ref{fastPricing} to solve pricing but solve the MP by adding one column at a time to $\Omega_R$ and solving RMP $\Psi(\Omega_R)$ to generate a dual solution.  We considered the Solomon instances data set \citep{solomon1987algorithms}, which consists of nine problem instances of various degrees of difficulty with 25 customers.

We compare the following parameterizations: a choice of either GM or standard CG with choices number of (0,4,6,8,10) LA-neighbors (for a total of 10 possibilities).  When using standard CG with 0 or 4 LA-neighbors the problem instances (3,4) and 4 respectively, can not be run to completion given a maximum run time of 24 hours.  

\begin{figure}
    \centering
\includegraphics[width=0.49\linewidth]{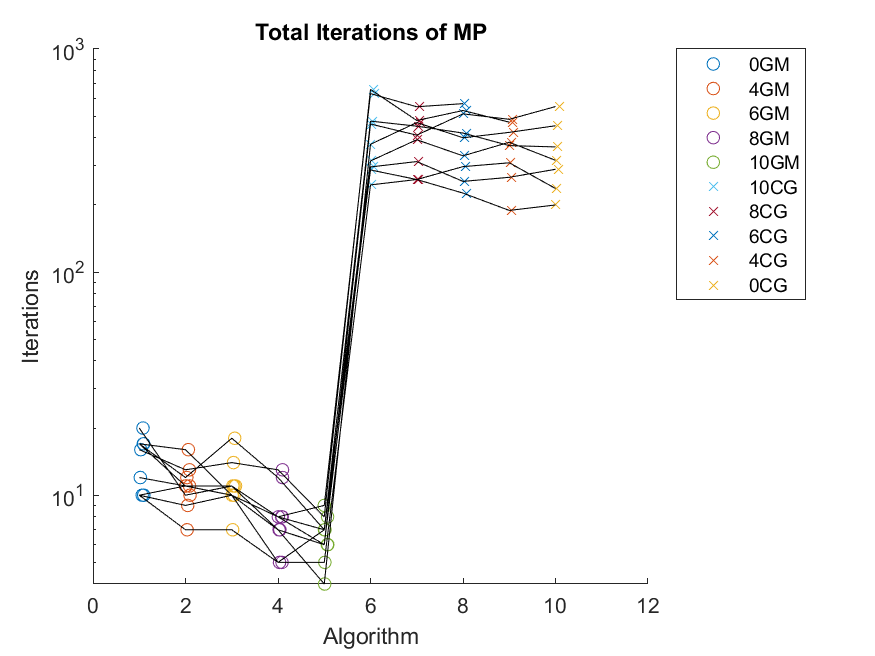} 
\includegraphics[width=0.49\linewidth]{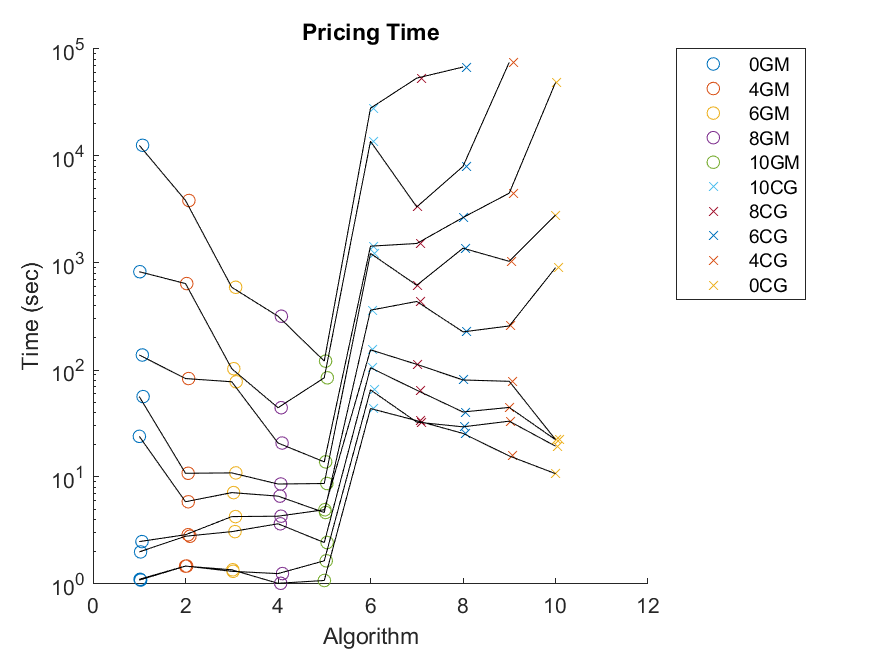}
\includegraphics[width=0.49\linewidth]{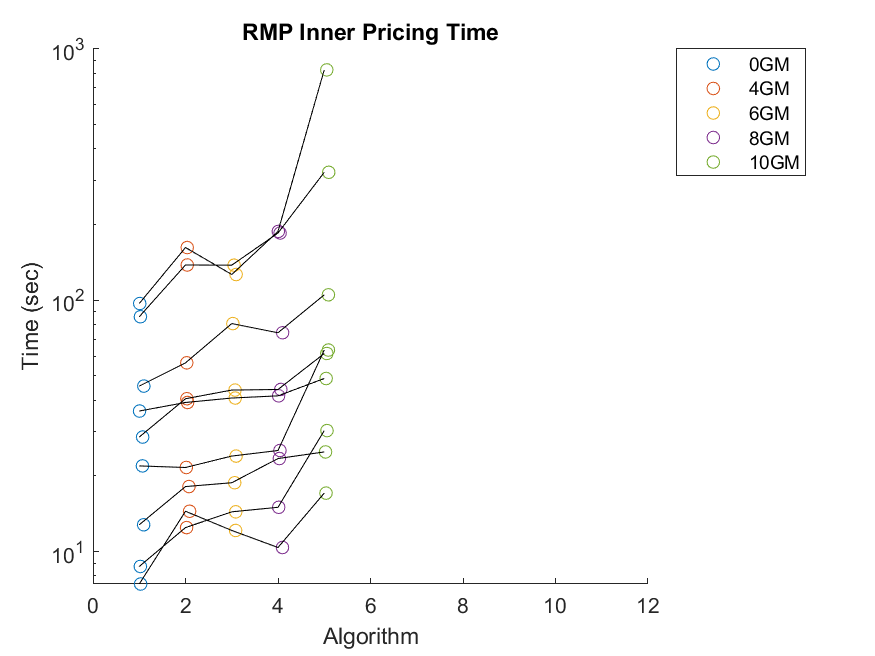}
\includegraphics[width=0.49\linewidth]{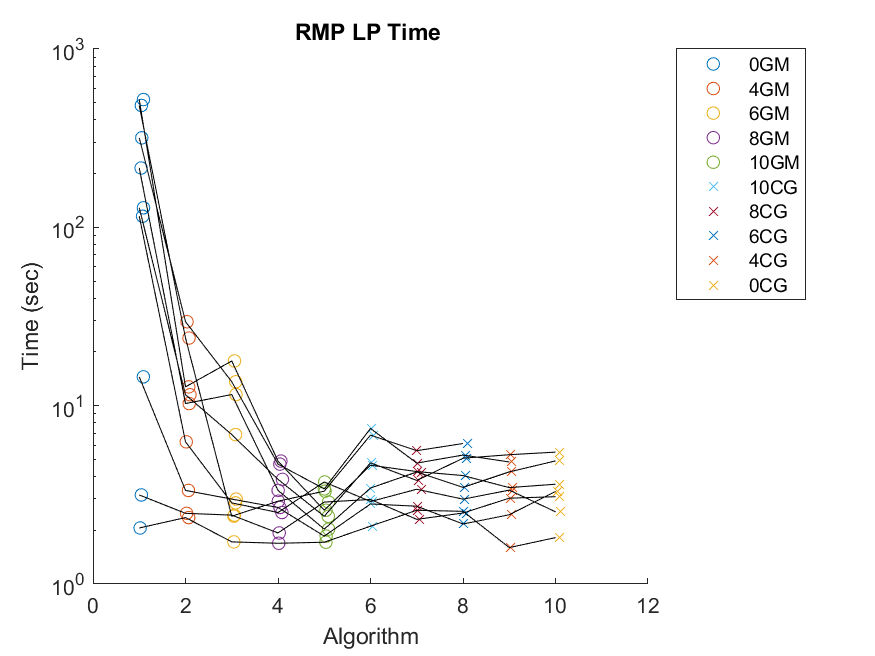}
\includegraphics[width=0.49\linewidth]{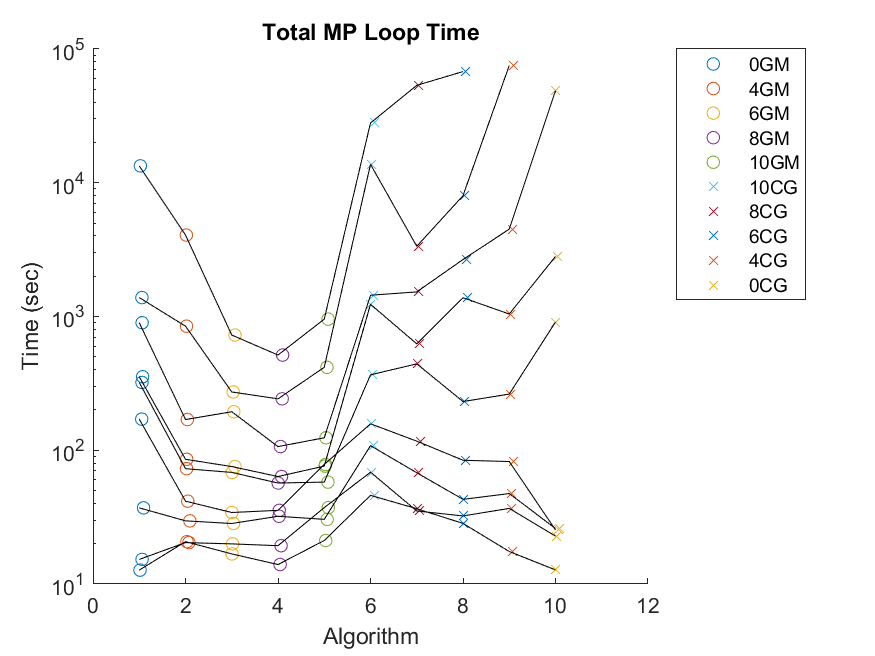}
\includegraphics[width=0.49\linewidth]{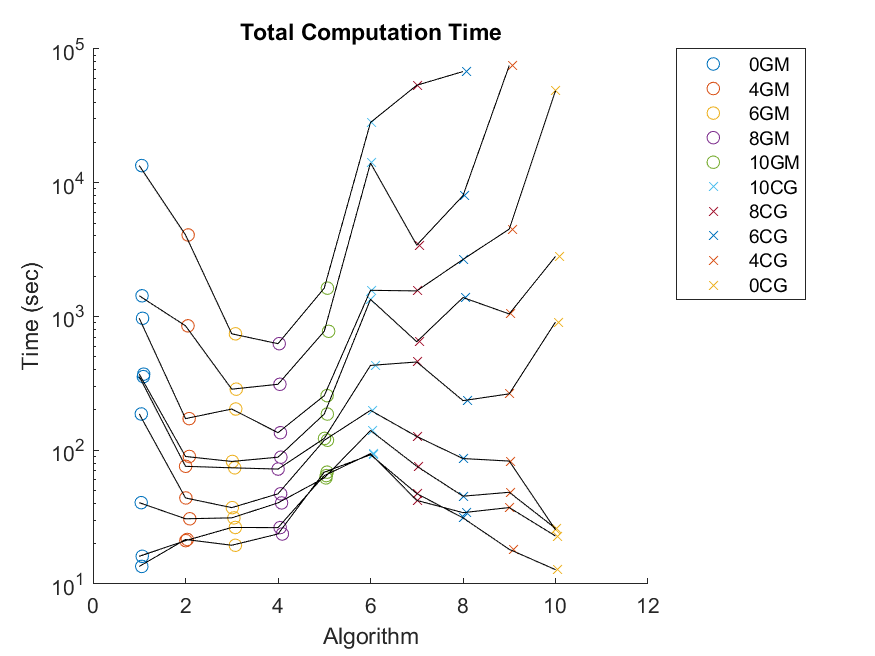}
    \caption{Performance of GM/CG under various numbers of LA-neighbors}
    \label{fig:my_results}
\end{figure}

All computation was done using MATLAB with default options.  Parallelization can be exploited in a future version of the code that was not done in our version.  Two key places to exploit parallelization are below.  We observe that most computation time during pricing is taken by computing $\hat{c}_{ij}$ in \eqref{hatCij}.  Since this computation can be easily parallelized we see room for further improvements in speed. Furthermore Line \ref{call_ez_pricing} which computes pricing over families is time intensive and can be sped up by solving Line \ref{call_ez_pricing} in parallel for each $l \in \Omega_R$.   We initialize $\Omega_R$ in standard CG with one column for each customer that services that single customer.  
In Fig \ref{fig:my_results} we plot the performance across our data set as a function of parameterization for the following performance metrics:  
\begin{itemize}
\item \textbf{Total Iterations of MP:}  This displays the number of iterations of the outer of GM (Lines \ref{loop_2}-\ref{loop_2_e} of Alg \ref{comp_alg}) or standard CG.  
\item \textbf{Pricing Time:}  This displays the time spent solving pricing over $\Omega$ over the course of Alg \ref{comp_alg}(Line \ref{call_fancy_pricing} of Alg \ref{comp_alg}) or standard CG.
\item \textbf{RMP Inner Pricing Time:}  This displays the time spent solving pricing over all $\Omega_l$ terms over the course of Alg \ref{comp_alg}(Line \ref{call_ez_pricing} of Alg \ref{comp_alg}). This step is not done in standard CG so the data is not plotted for standard CG.
\item \textbf{RMP LP Time:}  Total time spent solving the RMP of $\Psi(\Omega_R,I^-,P^- )$  in Alg \ref{comp_alg} ( Line  \ref{step_RMP} of Alg \ref{comp_alg} )  or  the RMP in standard CG.
\item \textbf{Total MP Loop Time:}  Total time spent in Alg \ref{comp_alg} or total time spent to solve the MP using standard CG. 
\item \textbf{Total Computation Time:  }Total time spent in Alg \ref{comp_alg} (or total time spent to solve the MP using standard CG) plus the time including pre-processing and post -processing.  Pre-processing corresponds to generating  $R^*_p$ via Alg \ref{better_ver}.   Post-processing corresponds to solving $\Psi(\Omega_R,I^-,P^- )$ as an ILP in order to generate a feasible integer solution.  All problem instances in our data set had tight LP relaxations.  
\end{itemize}
For each plot we use the x axis to describe parameterization.  The parameterization is characterized by the number of LA-neighbors followed by GM/or CG.  Note that $0$ LA-neighbors indicates that LA-arcs are not used.  We  plot one data point for each problem instance, parameterization.  We connect data points associated with a common problem instance with a black line.  

We observe large improvements with regards to  pricing time, total computation time, and number of iterations GM (Alg \ref{comp_alg}) when using GM over CG, and when using more LA-neighbors (up to 8 LA-neighbors with GM and 10 LA-neighbors with standard CG).  This is particularly true on the more time intensive problem instances.  Furthermore we observe that the improvements resulting from using GM and more LA-neighbors are complementary.  
\section{Conclusions and Future Work}
\label{sec_conc}
In this document we adapt Local Area (LA) routes \citep{mandal2022local_2} and Graph Master (GM)(\citep{yarkony2021graph,yarkony2022principled}) so as to permit the efficient solution to capacitated vehicle routing with time windows (CVRPTW) using column generation (CG).  Our GM approach projects each column generated, during the pricing step of CG, to a multi-graph where each path from source to sink corresponds to a feasible column.  Solving optimization over this restricted master problem (RMP) is efficient.  GM converges in far fewer iterations relative to standard CG.  Our approach for pricing adapts time window discretization \citep{boland2017continuous}, Decremental state space relaxation \citep{righini2009decremental} and LA-routes jointly.  Our pricing approach constructs a pricing graph upon which we iteratively tighten a relaxation of pricing till the shortest path corresponds to the lowest reduced cost column.  Our work can be generalized to other problems specifically problems where pricing is an (elementary) resource constrained shortest path problem as is common in operations research.  

In future work we intend to adapt LA-subset row inequalities (LA-SRI) \citep{mandal2022local_2} to efficiently tighten the underlying LP relaxation.  LA-SRI marginally weaken standard SRI \citep{jepsen2008subset} in such a manner that unlike original SRI allows for their efficient inclusion in pricing.  
We will also seek to incorporate our approach into a branch-cut-price \citep{barnprice} formulation.  In this case the standard branching rule of  \citep{ryan1981integer} can be used.  This  branches on the which customer/depot a given customer is succeeded by.  This fortunately can remove LA-arcs from consideration, which may speed convergence of optimization.  


\bibliographystyle{abbrvnat} 
\bibliography{col_gen_bib}

\begin{thebibliography}{26}
\providecommand{\natexlab}[1]{#1}
\providecommand{\url}[1]{\texttt{#1}}
\expandafter\ifx\csname urlstyle\endcsname\relax
  \providecommand{\doi}[1]{doi: #1}\else
  \providecommand{\doi}{doi: \begingroup \urlstyle{rm}\Url}\fi

\bibitem[Archetti et~al.(2011)Archetti, Bianchessi, and
  Speranza]{archetti2011column}
C.~Archetti, N.~Bianchessi, and M.~G. Speranza.
\newblock A column generation approach for the split delivery vehicle routing
  problem.
\newblock \emph{Networks}, 58\penalty0 (4):\penalty0 241--254, 2011.

\bibitem[Baldacci et~al.(2011)Baldacci, Mingozzi, and Roberti]{baldacci2011new}
R.~Baldacci, A.~Mingozzi, and R.~Roberti.
\newblock New route relaxation and pricing strategies for the vehicle routing
  problem.
\newblock \emph{Operations Research}, 59\penalty0 (5):\penalty0 1269--1283,
  2011.

\bibitem[Barnhart et~al.(1996)Barnhart, Johnson, Nemhauser, Savelsbergh, and
  Vance]{barnprice}
C.~Barnhart, E.~L. Johnson, G.~L. Nemhauser, M.~W.~P. Savelsbergh, and P.~H.
  Vance.
\newblock Branch-and-price: Column generation for solving huge integer
  programs.
\newblock \emph{Operations Research}, 46:\penalty0 316--329, 1996.

\bibitem[Boland et~al.(2017)Boland, Hewitt, Marshall, and
  Savelsbergh]{boland2017continuous}
N.~Boland, M.~Hewitt, L.~Marshall, and M.~Savelsbergh.
\newblock The continuous-time service network design problem.
\newblock \emph{Operations Research}, 65\penalty0 (5):\penalty0 1303--1321,
  2017.

\bibitem[Castro(2007)]{castro2007interior}
J.~Castro.
\newblock An interior-point approach for primal block-angular problems.
\newblock \emph{Computational optimization and Applications}, 36\penalty0
  (2-3):\penalty0 195--219, 2007.

\bibitem[Costa et~al.(2019)Costa, Contardo, and Desaulniers]{costa2019}
L.~Costa, C.~Contardo, and G.~Desaulniers.
\newblock Exact branch-price-and-cut algorithms for vehicle routing.
\newblock \emph{Transportation Science}, 26(1), 2019.

\bibitem[Desaulniers et~al.(2005)Desaulniers, Desrosiers, and
  Solomon]{Desaulniers2005}
G.~Desaulniers, J.~Desrosiers, and M.~M. Solomon, editors.
\newblock \emph{Column Generation}.
\newblock Springer, New York, 1st edition, 2005.

\bibitem[Desrochers et~al.(1992)Desrochers, Desrosiers, and
  Solomon]{Desrochers1992}
M.~Desrochers, J.~Desrosiers, and M.~Solomon.
\newblock A new optimization algorithm for the vehicle routing problem with
  time windows.
\newblock \emph{Operations Research}, 40\penalty0 (2):\penalty0 342--354, 1992.

\bibitem[Desrosiers and L{\"u}bbecke(2005)]{desrosiers2005primer}
J.~Desrosiers and M.~E. L{\"u}bbecke.
\newblock A primer in column generation.
\newblock In G.~Desaulniers, J.~Desrosiers, and M.~M. Solomon, editors,
  \emph{Column Generation}, pages 1--32. Springer, New York, NY, 2005.

\bibitem[Du~Merle et~al.(1999)Du~Merle, Villeneuve, Desrosiers, and
  Hansen]{du1999stabilized}
O.~Du~Merle, D.~Villeneuve, J.~Desrosiers, and P.~Hansen.
\newblock Stabilized column generation.
\newblock \emph{Discrete Mathematics}, 194\penalty0 (1-3):\penalty0 229--237,
  1999.

\bibitem[Geoffrion(1974)]{geoffrion1974lagrangean}
A.~M. Geoffrion.
\newblock Lagrangean relaxation for integer programming.
\newblock In \emph{Approaches to integer programming}, pages 82--114. Springer,
  1974.

\bibitem[Gilmore and Gomory(1961)]{cuttingstock}
P.~Gilmore and R.~Gomory.
\newblock A linear programming approach to the cutting-stock problem.
\newblock \emph{Operations Research}, 9\penalty0 (6):\penalty0 849--859, 1961.

\bibitem[Haghani et~al.(2021)Haghani, Li, Koenig, Kunapuli, Contardo, Regan,
  and Yarkony]{haghani2021multi}
N.~Haghani, J.~Li, S.~Koenig, G.~Kunapuli, C.~Contardo, A.~Regan, and
  J.~Yarkony.
\newblock Multi-robot routing with time windows: A column generation approach.
\newblock \emph{arXiv preprint arXiv:2103.08835}, 2021.

\bibitem[Irnich and Desaulniers(2005)]{irnich2005shortest}
S.~Irnich and G.~Desaulniers.
\newblock Shortest path problems with resource constraints.
\newblock In G.~Desaulniers, J.~Desrosiers, and M.~M. Solomon, editors,
  \emph{Column generation}, pages 33--65. Springer, 2005.

\bibitem[Jepsen et~al.(2008)Jepsen, Petersen, Spoorendonk, and
  Pisinger]{jepsen2008subset}
M.~Jepsen, B.~Petersen, S.~Spoorendonk, and D.~Pisinger.
\newblock Subset-row inequalities applied to the vehicle-routing problem with
  time windows.
\newblock \emph{Operations Research}, 56\penalty0 (2):\penalty0 497--511, 2008.

\bibitem[Mandal et~al.(2022)Mandal, Regan, and Yarkony]{mandal2022local_2}
U.~Mandal, A.~Regan, and J.~Yarkony.
\newblock Local area routes and valid inequalities for efficient vehicle
  routing.
\newblock \emph{arXiv preprint arXiv:2209.12963}, 2022.

\bibitem[Marsten et~al.(1975)Marsten, Hogan, and
  Blankenship]{marsten1975boxstep}
R.~E. Marsten, W.~Hogan, and J.~W. Blankenship.
\newblock The boxstep method for large-scale optimization.
\newblock \emph{Operations Research}, 23\penalty0 (3):\penalty0 389--405, 1975.

\bibitem[Pessoa et~al.(2018)Pessoa, Sadykov, Uchoa, and
  Vanderbeck]{Pessoa2018Automation}
A.~A. Pessoa, R.~Sadykov, E.~Uchoa, and F.~Vanderbeck.
\newblock Automation and combination of linear-programming based stabilization
  techniques in column generation.
\newblock \emph{{INFORMS} Journal on Computing}, 30\penalty0 (2):\penalty0
  339--360, 2018.
\newblock \doi{10.1287/ijoc.2017.0784}.

\bibitem[Righini and Salani(2008)]{righini2008new}
G.~Righini and M.~Salani.
\newblock New dynamic programming algorithms for the resource constrained
  elementary shortest path problem.
\newblock \emph{Networks: An International Journal}, 51\penalty0 (3):\penalty0
  155--170, 2008.

\bibitem[Righini and Salani(2009)]{righini2009decremental}
G.~Righini and M.~Salani.
\newblock Decremental state space relaxation strategies and initialization
  heuristics for solving the orienteering problem with time windows with
  dynamic programming.
\newblock \emph{Computers \& Operations Research}, 36\penalty0 (4):\penalty0
  1191--1203, 2009.

\bibitem[Ropke and Cordeau(2009)]{ropke2009branch}
S.~Ropke and J.-F. Cordeau.
\newblock Branch and cut and price for the pickup and delivery problem with
  time windows.
\newblock \emph{Transportation Science}, 43\penalty0 (3):\penalty0 267--286,
  2009.

\bibitem[Rousseau et~al.(2007)Rousseau, Gendreau, and
  Feillet]{rousseau2007interior}
L.-M. Rousseau, M.~Gendreau, and D.~Feillet.
\newblock Interior point stabilization for column generation.
\newblock \emph{Operations Research Letters}, 35\penalty0 (5):\penalty0
  660--668, 2007.

\bibitem[Ryan and Foster(1981)]{ryan1981integer}
D.~M. Ryan and B.~A. Foster.
\newblock An integer programming approach to scheduling.
\newblock \emph{Computer scheduling of public transport urban passenger vehicle
  and crew scheduling}, pages 269--280, 1981.

\bibitem[Solomon(1987)]{solomon1987algorithms}
M.~M. Solomon.
\newblock Algorithms for the vehicle routing and scheduling problems with time
  window constraints.
\newblock \emph{Operations research}, 35\penalty0 (2):\penalty0 254--265, 1987.

\bibitem[Yarkony and Regan(2022)]{yarkony2022principled}
J.~Yarkony and A.~Regan.
\newblock Principled graph management.
\newblock \emph{arXiv preprint arXiv:2202.01274}, 2022.

\bibitem[Yarkony et~al.(2021)Yarkony, Haghani, and Regan]{yarkony2021graph}
J.~Yarkony, N.~Haghani, and A.~Regan.
\newblock Graph generation: A new approach to solving expanded linear
  programming relaxations.
\newblock \emph{arXiv preprint arXiv:2110.01070}, 2021.

\end{thebibliography}

\appendix

\section{Proof of Equivalent Representation of Departure Time}
\label{sec_proof}

In this section we prove that \eqref{tauDef} accurately characterizes $T_{r}(t)$ for all $r\in R_p,t_0\geq t \geq 0, p \in P$.  We prove this using induction. Observe that for the base case where $|N_p|=0$ that the claim holds by definition.  If \eqref{tauDef} is does not hold for all cases then there must exist some $r,t$ be defined so that \eqref{tauDef} but \eqref{tauDef} holds for $r^-$ and all $t_0 \geq t\geq 0$.  We describe this formally below using $u,w$ to denote the first two customers in $r$. 

\textbf{Claim to be Proven False}:  
\begin{align}
\label{claim_to_false}
T_{r}(t)=T_{r^-}(-t_{uw}+\min(t,t^+_u))-\infty*[t<t^-_u] \neq -c_r+\min(t,\tau_{r1}) -\infty[\min(t,t^+_u)<\tau_{r2}]
\end{align}
For now we assume \eqref{claim_to_false} is true.  

\textbf{Proof:  }
We now use \eqref{tauDef} to replace $T_{r^-}(-t_{uw}+\min(t,t^+_u))$ in \eqref{claim_to_false}.  
\begin{align}
\label{eq_term_1}
c_{r^-}+\min(-t_{uw}+\min(t,t^+_u),\tau_{r^{-}1}) -\infty[-t_{uw}+\min(t,t^+_u)<\tau_{r^-2}]\\
\neq  -c_r+\min(t,\tau_{r1}) -\infty[\min(t,t^+_u)<\tau_{r2}]\nonumber 
\end{align}
We now consider the terms corresponding to $\infty$ and those not corresponding to $\infty$ separately in \eqref{eq_term_1}. For the \eqref{claim_to_false} to be true then it must be the case that  $T_{r^-}(-t_{uw}+\min(t,t^+_u))-\infty*[t<t^-_u] \neq -c_r+\min(t,\tau_{r1}) -\infty[\min(t,t^+_u)<\tau_{r2}]$.  A necessary condition for \eqref{claim_to_false} to be true is that  one or both of the following must must be satisfied.  
\begin{subequations}
\label{2Pos}
\begin{align}
-\infty[\min(t,t^+_u)<\tau_{r2}]\neq -\infty[-t_{uw}+\min(t,t^+_u)<\tau_{r^-2}]-\infty*[t<t^-_u]\label{pos_1a}\\
 -c_{r^-} +\min(-t_{uw}+\min(t,t^+_u),\tau_{r^{-}1}) -\neq c_r + \min(t,\tau_{r1}) \label{pos_2a}
\end{align}
\end{subequations}

Let us consider \eqref{pos_1a} first.  We apply the following transformations including plugging in the definition of $\tau_{r2}$.  
\begin{subequations}
\begin{align}
\label{infEq}
-\infty[\min(t,t^+_u)<\tau_{r2}]\neq -\infty[-t_{uw}+\min(t,t^+_u)<\tau_{r^-2}]-\infty*[t<t^-_u]\\
-\infty[\min(t,t^+_u)<\max(t^-_u,\tau_{r^-2}+t_{uw})]\neq -\infty[-t_{uw}+\min(t,t^+_u)<\tau_{r^-2}]-\infty*[t<t^-_u]\\
-\infty[\min(t,t^+_u)<\max(t^-_u,\tau_{r^-2}+t_{uw})]\neq -\infty[\min(t,t^+_u)<\tau_{r^-2}+t_{uw}]-\infty*[t<t^-_u]\\
-\infty[\min(t,t^+_u)<\max(t^-_u,\tau_{r^-2}+t_{uw})]\neq -\infty[\min(t,t^+_u)<\tau_{r^-2}+t_{uw}]-\infty*[t<t^-_u] \label{finEqAset}
\end{align}
\end{subequations}
Consider that  $t<t^-_u$  then both the LHS and RHS of \eqref{finEqAset} equal $-\infty$.  Thus if a violation occurs in \eqref{finEqAset} then $t\geq t^-_u$.  Thus we gain the following expression for a potential violation:  
\begin{align}
-\infty[\min(t,t^+_u)<\max(t^-_u,\tau_{r^-2}+t_{uw})]\neq -\infty[\min(t,t^+_u)<\tau_{r^-2}+t_{uw}] \label{tryme}
\end{align}
The only way for \eqref{tryme} to produce an inequality is if $\tau_{r^-2}+t_{uw}\leq \min(t,t^+_u)<t^-_u$.  However both $t$ and $t^+_u$ are greater than or equal to $t^-_u$ creating a contradiction.  
However the two LHS and RHS of \eqref{infEq} are identical so \eqref{pos_1a} does not hold.  

We now consider the case where \eqref{pos_2a} is violated.  We now plug in the definition of $\tau_{r1}$, into \eqref{pos_2a}, and  noting that $t_{uw}=c_{uw}$ to produce the following expressions by equivalent transformations.  
\begin{subequations}
\begin{align}
-c_{r^-}+\min(-t_{uw}+\min(t,t^+_u),\tau_{r^{-}1})\neq -c_r+\min(t,\tau_{r1})\\
 -c_{r^-}+\min(-t_{uw}+\min(t,t^+_u),\tau_{r^{-}1})\neq -c_{r^-}-c_{uw}+\min(t,\min(t^+_u,\tau_{r^-1}+t_{uw}))\\
 \min(\min(t-t_{uw},t^+_u-t_{uw}),\tau_{r^{-}1})\neq -c_{uw}+\min(t,\min(t^+_u,\tau_{r^-1}+t_{uw}))\\
 \min(\min(t-t_{uw},t^+_u-t_{uw}),\tau_{r^{-}1})\neq -t_{uw}+\min(t,\min(t^+_u,\tau_{r^-1}+t_{uw}))\\
\min(t-t_{uw},t^+_u-t_{uw},\tau_{r^{-}1})\neq -t_{uw}+\min(t,t^+_u,\tau_{r^-1}+t_{uw})\\
 \min(t-t_{uw},t^+_u-t_{uw},\tau_{r^{-}1}) \neq \min(t-t_{uw},t^+_u-t_{uw},\tau_{r^-1})\label{finverGroup}
\end{align}
\end{subequations}

Since the RHS and LHS of \eqref{finverGroup} are identical so the inequality does not hold \eqref{pos_2a} does not hold.  Since neither \eqref{pos_1a} or \eqref{pos_2a} hold then \eqref{claim_to_false} is false.  Thus \eqref{tauDef} must be true.

\end{document}